\documentclass{amsart}
\newtheorem{thm}{Theorem}[section]
\newtheorem{cor}[thm]{Corollary}
\newtheorem{lem}[thm]{Lemma}

\newtheorem{claim}[thm]{Claim}
\newtheorem{prop}[thm]{Proposition}

\newtheorem{rem}[thm]{Remark}
\newtheorem{ex}[thm]{Example}

\newtheorem{op}[thm]{Open Problem}
\newtheorem{defn}[thm]{Definition}
\newtheorem{notation}[thm]{Notation}
\numberwithin{equation}{section}

\begin{document}

\title[Nonsmooth analysis on Riemannian manifolds]{Nonsmooth analysis
and Hamilton-Jacobi equations on Riemannian manifolds}
\author{Daniel Azagra, Juan Ferrera, Fernando L{\'o}pez-Mesas.}

\address{Departamento de An{\'a}lisis Matem{\'a}tico\\ Facultad de
Matem{\'a}ticas\\ Universidad Complutense\\ 28040 Madrid, Spain}

\date{May 23, 2003}

\email{Daniel\_Azagra@mat.ucm.es; ferrera@mat.ucm.es; FLopez\_Mesas@mat.ucm.es}

\thanks{The first-named author was supported by a
Marie Curie Fellowship of the European Community Training and
Mobility of Researchers Programme under contract number
HPMF-CT-2001-01175. The second-named author was partially
supported by BFM grant 2000/0609. The third-named author enjoyed
a study license granted by the Consejer{\'\i}a de Educaci{\'o}n
de la Comunidad de Madrid.}



\begin{abstract}
We establish some perturbed minimization principles, and we
develop a theory of subdifferential calculus, for functions
defined on Riemannian manifolds. Then we apply these results to
show existence and uniqueness of viscosity solutions to
Hamilton-Jacobi equations defined on Riemannian manifolds.
\end{abstract}

\maketitle

\section{Introduction}

The aim of this paper is threefold. First, we extend some
perturbed minimization results such as the smooth variational
principle of Deville, Godefroy and Zizler, and other
almost-critical-point-spotting results, such as approximate
Rolle's type theorems, to the realm of Riemannian manifolds.
Second, we introduce a definition of subdifferential for functions
defined on Riemannian manifolds, and we develop a theory of
subdifferentiable calculus on such manifolds that allows most of
the known applications of subdifferentiability to be extended to
Riemannian manifolds. For instance, we show that every convex
function on a Riemannian manifold (that is, every function which
is convex along geodesics) is everywhere subdifferentiable (on the
other hand, every continuous function is superdifferentiable on a
dense set, hence convex functions are differentiable on dense
subsets of their domains). Third, we also use this theory to prove
existence and uniqueness of viscosity solutions to Hamilton-Jacobi
equations defined on Riemannian manifolds. Let us introduce some
of these results.

It is known that the classic Rolle's theorem fails in
infinite-dimensions, that is, in every infinite-dimensional
Banach space with a $C^1$ smooth (Lipschitz) bump function there
are $C^1$ smooth (Lipschitz) functions which vanish outside a
bounded open set and yet have a nonzero derivative everywhere
inside this set; see \cite{AJ1} and the references therein. In
fact, the failure of Rolle's theorem infinite dimensions takes on
a much more dramatic form in a recent result of D. Azagra and M.
Cepedello Boiso \cite{AC2}: the smooth functions {\em with no
critical points} are dense in the space of continuous functions
on every Hilbert manifold (this result may in turn be viewed as a
very strong approximate version for infinite dimensional
manifolds of the Morse-Sard theorem). So, when we are given a
smooth function on an infinite-dimensional Riemannian manifold we
should not expect to be able to find any critical point, whatever
the overall shape of this function is, as there might be none.
This important difference between finite and infinite dimensions
forces us to consider approximate substitutes of Rolle's theorem
and the classic minimization principles, looking for the
existence of arbitrarily small derivatives (instead of vanishing
ones) for every function satisfying (in an approximate manner)
the conditions of the classical Rolle's theorem. This is what the
papers \cite{AGJ, AD1} deal with, one in the differentiable case
(showing for instance that if a differentiable function
oscillates less than $2\varepsilon$ on the boundary of a unit
ball then there is a point inside the ball such that the
derivative of the function has norm less than or equal to
$\varepsilon$), and the other in the subdifferentiable one. More
generally, a lot of perturbed minimization (or variational)
principles have been studied, perhaps the most remarkable being
Ekeland's Principle, Borwein-Preiss' Principle, and
Deville-Godefroy-Zizler's Smooth Variational Principle. See
\cite{DGZ, DeGou} and the references therein.

There are many important applications of those variational
principles. Therefore, it seems reasonable to look for analogues
of these perturbed minimizations principles within the theory of
Riemannian manifolds. In Section 3 we prove some
almost-critical-point spotting results. First we establish an
approximate version of Rolle's theorem which holds for
differentiable mappings defined on subsets of arbitrary Riemannian
manifolds. Then we give a version of Deville-Godefroy-Zizler
smooth variational principle which holds for those complete
Riemannian manifold $M$ which are {\em uniformly bumpable}
(meaning that there exists some numbers $R>1, r>0$ such that for
every point $p\in M$ and every $\delta\in (0,r)$ there exists a
function $b:M\to [0,1]$ such that $b(x)=0$ if $d(x,p)\geq\delta$,
$b(p)=1$, and $\sup_{x\in M}\|db(x)\|_{x}\leq R/\delta$. Of course
every Hilbert space is uniformly bumpable, and there are many
other examples of uniformly bumpable manifolds: as we will see,
every Riemannian manifold which has strictly positive injectivity
and convexity radii is uniformly bumpable). For those Riemannian
manifolds we show that, for every lower semicontinuous function
$f:M\to (-\infty, \infty]$ which is bounded below, there exists a
$C^1$ smooth function $\varphi:M\to\mathbb{R}$, which is
arbitrarily small and has an arbitrarily small derivative
everywhere, such that $f-\varphi$ attains a strong global minimum
at some $p\in M$.

This result leads up to one of the main topic of this paper:
subdifferentiability of functions on Riemannian manifolds, since,
according to the definition we are going to give of
subdifferential, this implies that such $f$ is subdifferentiable
at the point $p$. We will say that a function $f:M\to (-\infty,
\infty]$ is subdifferentiable at $p$ provided there exists a $C^1$
smooth function $\varphi:M\to\mathbb{R}$ such that $f-\varphi$
attains a local minimum at $p$. The set of the derivatives
$d\varphi(p)$ of all such functions $\varphi$ will be called
subdifferential of $f$ at $p$, a subset of $T^{*}M_{p}$ which will
be denoted by $D^{-}f(p)$. Of course, when $M$ is $\mathbb{R}^{n}$
or a Hilbert space, this definition agrees with the usual one.
Apart from being a useful generalization of the theory of
subdifferentiability of convex functions, this notion of
subdifferentiability plays a fundamental role in the study of
Hamilton-Jacobi equations in $\mathbb{R}^{n}$ and
infinite-dimensional Banach spaces. Not only is this concept
necessary to understand the notion of {\em viscosity solution}
(introduced by M. G. Crandall and P. L. Lions, see \cite{CL1, CL2,
CL3, CL4, CL5, CL6, CL7, CL8, CL9, CL10}); from many results
concerning subdifferentials one can also deduce relatively easy
proofs of the existence, uniqueness and regularity of viscosity
solutions to Hamilton-Jacobi equations; see, for instance,
\cite{DeTVM, De2, DH, DH2}. We refer to \cite{Devillesolosurvey,
DeGou} for an introduction to subdifferential calculus in Banach
spaces and its applications (especially Hamilton-Jacobi
equations).

Section 4 is devoted to the study of subdifferentials of functions
defined on manifolds. We start by giving other equivalent
definitions of subdifferentiability and superdifferentiability,
including a local one through charts, which sometimes makes it
easy to translate some results already established in the
$\mathbb{R}^{n}$ or the Banach space cases to the setting of
Riemannian manifolds. We also show that a function $f$ is
differentiable at a point $p$ if and only if $f$ is both
subdifferentiable and superdifferentiable at $p$. Next we study
the elementary properties of this subdifferential with respect to
sums, products and composition, including direct and inverse fuzzy
rules. We finish this section by establishing two mean value
theorems, and showing that lower semicontinuous functions are
subdifferentiable on dense subsets of their domains.

In section 5 we study the links between convexity and
(sub)differentiability of functions defined on Riemannian
manifolds. Recall that a function $f:M\to\mathbb{R}$ defined on a
Riemannian manifold $M$ is said to be convex provided
$f\circ\sigma$ is convex, for every geodesic $\sigma$. The papers
\cite{GreeneShio1, GreeneShio2, GreeneWu1, GreeneWu2} provide a
very good introduction to convexity on Riemannian manifolds and
the geometrical implications of the existence of global convex
functions on a Riemannian manifold; for instance it is shown in
\cite{GreeneShio1} that every two-dimensional manifold which
admits a global convex function which is locally nonconstant must
be diffeomorphic to the plane, the cylinder, or the open
M{\"o}bius strip. Among other things, we show in this section
that every convex function defined on a Riemannian manifold is
everywhere subdifferentiable, and is differentiable on a dense
set (when the manifold is finite-dimensional, the set of points
of nondifferentiability has measure zero).

Finally, in Section 6 we study some Hamilton-Jacobi equations
defined on Riemannian manifolds (either finite or
infinite-dimensional). Examples of Hamilton-Jacobi equations
arise naturally in the setting of Riemannian manifolds, see
\cite{AbrahamMarsden} in relation to Lyapounov theory and optimal
control. However, we do not know of any work that has studied
nonsmooth solutions, in general, or viscosity solutions, in
particular, to Hamilton-Jacobi equations defined on Riemannian
manifolds . This may be due to the lack of a theory of nonsmooth
calculus for functions defined on Riemannian manifolds. Here we
will show how the subdifferential calculus and perturbed
minimization principles that we develop in the previous sections
can be applied to get results on existence and uniqueness of
viscosity solutions to equations of the form $$\left\{
\begin{matrix}u+F(du)=f \cr u \textrm{ bounded, } \cr\end{matrix} \right.$$
where $f:M\to \mathbb{R}$ is a bounded uniformly continuous
function, and $F:T^*M\to \mathbb{R}$ is a function defined on the
cotangent bundle of $M$ which satisfies a uniform continuity
condition. The manifold $M$ must also satisfy that it has positive
convexity and injectivity radii (this condition is automatically
met by every compact manifold, for instance). We also prove some
results about ``regularity" (meaning differentiability almost
everywhere) of the viscosity solutions to some of these equations.
Finally, we study the equation $\|du(x)\|_{x}=1$ for all
$x\in\Omega$, $u(x)=0$ for all $x\in\partial\Omega$, where
$\Omega$ is a bounded open subset of $M$, and we show that
$x\mapsto d(x,\partial\Omega)$ is the unique viscosity solution to
this equation (which has no classical solution).

Let us also mention that one can develop other theories of
subdifferentiability on Riemannian manifolds which have very
interesting applications as well. For instance, one can consider
quite natural extensions of the proximal subgradient and Clarke's
generalized gradient for functions defined on Riemannian
manifolds. By teaming these notions with infimal convolutions with
squared geodesic distance functions one can get very interesting
results about differentiability properties of the distance
function to a closed subset of a Riemannian manifold; see
\cite{AFL3}.

\medskip


\section{Preliminaries and tools}

In this section we recall some definitions and known results about
Riemannian manifolds which will be used later on.

We will be dealing with functions defined on Riemannian manifolds
(either finite or infinite-dimensional). A Riemannian manifold
$(M,g)$ is a $C^\infty$ smooth manifold $M$ modelled on some
Hilbert space $H$ (possibly infinite-dimensional), such that for
every $p\in M$ we are given a scalar product
$g(p)=g_{p}:=\langle\cdot,\cdot\rangle_{p}$ on the tangent space
$TM_{p}\simeq H$ so that $\|x\|_{p}=(\langle
x,x\rangle_{p})^{1/2}$ defines an equivalent norm on $TM_{p}$ for
each $p\in M$, and in such a way that the mapping $p\in M \mapsto
g_{p}\in\mathcal{S}^{2}(M)$ is a $C^\infty$ section of the bundle
$\sigma_{2}:\mathcal{S}_{2}\to M$ of symmetric bilinear forms.

If a function $f:M\longrightarrow\mathbb{R}$ is differentiable at
$p\in M$, the norm of the differential $df(p)\in T^{*}M_{p}$ at
the point $p$ is defined by
    $$
    \|df(p)\|_{p}=\sup\{df(p)(v) : v\in TM_{p}, \|v\|_{p}\leq 1\}.
    $$
Since $(TM_{p}, \|\cdot\|_{p})$ is a Hilbert space, we have a
linear isometric identification between this space and its dual
$(T^{*}M_{p}, \|\cdot\|_{p})$ through the mapping $TM_{p}\ni
x\mapsto f_{x}=x\in T^{*}M_{p}$, where $f_{x}(y)=\langle x,
y\rangle$ for every $y\in TM_{p}$.

For every piecewise $C^1$ smooth path $\gamma:[a,b]\to M$ we
define its length as
    $$
    L(\gamma)=
    \int_{a}^{b}\|\frac{d\gamma}{dt}(s)\|_{\gamma(s)}ds.
    $$
This length depends only on the path $\gamma[a,b]$ itself, and
not on the way the point $\gamma(t)$ moves along it: if
$h:[0,1]\to[a,b]$ is a continuous monotone function then
$L(\gamma\circ h)=L(\gamma)$. We can always assume that a path
$\gamma$ is parameterized by arc length, which means that
$\gamma:[0, T]\to M$ satisfies
$\|\frac{d\gamma}{dt}(s)\|_{\gamma(s)}=1$ for all $s$, and
therefore
    $$
    L(\gamma_{|_{[0,r]}})=\int_{0}^{r}\|\frac{d\gamma}{dt}(s)\|_{\gamma(s)}ds=r
    $$
for each $r\in[0,T]$. For any two points $p,q\in M$, let us define
    $$
    d(p,q)=\inf\{ L(\gamma): \gamma \text{ is a } C^{1} \text{ smooth path joining }
    p \text{ and } q \text{ in } M\}.
    $$
Then $d$ is a metric on $M$ (called the $g$-distance on $M$) which
defines the same topology as the one $M$ naturally has as a
manifold. For this metric we define the closed ball of center $p$
and radius $r>0$ as
    $$
    B_{g}(p,r)=\{q\in M: d(p,q)\leq r\}.
    $$
Let us recall that in every Riemannian manifold there is a unique
natural covariant derivation, namely the Levi-Civita connection
(see Theorem 1.8.11 of \cite{Klingenberg}); following Klingenberg
we denote this derivation by $\nabla_{X} Y$ for any vector fields
$X, Y$ on $M$. We should also recall that a geodesic is a
$C^\infty$ smooth path $\gamma$ whose tangent is parallel along
the path $\gamma$, that is, $\gamma$ satisfies the equation
$\nabla_{d\gamma(t)/dt}d\gamma(t)/dt=0$. A geodesic always
minimizes the distance between points which are close enough to
each other.

Any path $\gamma$ joining $p$ and $q$ in $M$ such that
$L(\gamma)=d(p,q)$ is a geodesic, and it is called a minimal
geodesic. In the sequel all geodesic paths will be assumed to be
parameterized by arc length, unless otherwise stated.

\begin{thm}[Hopf-Rinow]
If $M$ is a finite-dimensional Riemannian manifold which is
complete and connected, then there is at least one minimal
geodesic connecting any two points in $M$.
\end{thm}
On the other hand, for any given point $p$, the statement ''$q$
can be joined to $p$ by a unique minimal geodesic`` holds for
almost every $q\in M$; see \cite{Milnor}.

As is well known, the Hopf-Rinow theorem fails when $M$ is
infinite-dimensional, but Ekeland \cite{Ekeland2} proved (by
using his celebrated variational principle) that, even in infinite
dimensions, the set of points that can be joined by a minimal
geodesic in $M$ is dense.

\begin{thm}[Ekeland]\label{Ekelan's approximated Hopf-Rinow theorem}
If $M$ is an infinite-dimensional Riemannian manifold which is
complete and connected then, for any given point $p$, the set
$\{q\in M : q \textrm{ {\em can be joined to $p$ by a unique
minimal geodesic}}\}$ is residual in $M$.
\end{thm}

\medskip

The existence theorem for ODEs implies that for every $V\in TM$
there is an open interval $J(V)$ containing $0$ and a unique
geodesic $\gamma_{V}:J(V)\to M$ with $d\gamma(0)/dt=V$. This in
turn implies that there is an open neighborhood $\tilde{T}M$ of
$M$ in $TM$ such that for every $V\in \tilde{T}M$, the geodesic
$\gamma_{V}(t)$ is defined for $|t|<2$. The exponential mapping
$\exp:\tilde{T}M\to M$ is then defined as $\exp
(V)=\gamma_{V}(1)$, and the restriction of $\exp$ to a fiber
$TM_{x}$ in $\tilde{T}M$ is denoted by $\exp_{x}$.

Let us now recall some useful properties of the exponential map.
See \cite{Lang, Klingenberg}, for instance, for a proof of the
following theorem.
\begin{thm}\label{properties of exp}
For every Riemannian manifold $(M,g)$ and every $x\in M$ there
exists a number $r>0$ and a map $\exp_{x}:B(0_{x},r)\subset
TM_{x}\to M$ such that
\begin{enumerate}
\item $\exp_{x}:B(0_{x},\delta)\to B_{M}(x,\delta)$ is a bi-Lipschitz
$C^\infty$ diffeomorphism, for all $\delta\in (0,r]$.
\item $\exp_{x}$ takes the segments passing through $0_{x}$ and
contained in $B(0_{x},r)\subset TM_{x}$ into geodesic paths in
$B_{M}(x,r)$.
\item $d\exp_{x}(0_{x})=\textrm{id}_{TM_{x}}$.
\end{enumerate}
In particular, taking into account condition $(3)$, for every
$C>1$, the radius $r$ can be chosen to be small enough so that the
mappings $\exp_{x}:B(0_{x},r)\to B_{M}(x,r)$ and
$\exp_{x}^{-1}:B_{M}(x, r)\to B(0_{x},r)$ are $C$-Lipschitz.
\end{thm}

Recall that a Riemannian manifold $M$ is said to be {\em
geodesically complete} provided the maximal interval of
definition of every geodesic in $M$ is all of $\mathbb{R}$. This
amounts to saying that for every $x\in M$, the exponential map
$\exp_{x}$ is defined on all of the tangent space $TM_{x}$
(though, of course, $\exp_{x}$ is not necessarily injective on all
of $TM_{x}$). It is well known that every complete Riemannian
manifold is geodesically complete. In fact we have the following
result (see \cite{Lang}, p. 224 for a proof).
\begin{prop}
Let $(M,g)$ be a Riemannian manifold. Consider the following
conditions:
\begin{enumerate}
\item $M$ is complete (with respect to the $g$-distance).
\item All geodesics in $M$ are defined on $\mathbb{R}$.
\item For every $x\in M$, the exponential map $\exp_{x}$ is
defined on all of $TM_{x}$.
\item There is some $x\in M$ such that the exponential map
$\exp_{x}$ is defined on all of $TM_{x}$.
\end{enumerate}
Then, $(1)\implies (2)\implies (3)\implies (4)$. Furthermore, if
we assume that $M$ is {\em finite-dimensional}, then all of the
four conditions are equivalent to a fifth:
\begin{enumerate}
\item[{(5)}] Every closed and $d_{g}$-bounded subset of $M$ is
compact.
\end{enumerate}
\end{prop}

\noindent Next let us recall some results about convexity in
Riemannian manifolds.

\begin{defn}
{\em We say that a subset $U$ of a Riemannian manifold is {\em
convex} if given $x, y\in U$ there exists a unique geodesic in
$U$ joining $x$ to $y$, and such that the length of the geodesic
is $\textrm{dist}(x,y)$.}
\end{defn}
Every Riemannian manifold is {\em locally convex}, in the
following sense.
\begin{thm}[Whitehead]\label{Riemannian manifolds are locally convex}
Let $M$ be a Riemannian manifold. For every $x\in M$, there
exists $c>0$ such that for all $r$ with $0<r<c$, the open ball
$B(x,r)=\exp_{x} B(0_{x}, r)$ is convex.
\end{thm}
This theorem gives rise to the notion of {\em uniformly locally
convex} manifold, which will be of interest when discussing smooth
variational principles and Hamilton-Jacobi equations on
Riemannian manifolds.
\begin{defn}\label{uniformly locally convex manifold}
{\em We say that a Riemannian manifold $M$ is {\em uniformly
locally convex} provided that there exists $c>0$ such that for
every $x\in M$ and every $r$ with $0<r<c$ the ball
$B(x,r)=\exp_{x} B(0_{x}, r)$ is convex.}
\end{defn}
This amounts to saying that the global convexity radius of $M$ (as
defined below) is strictly positive.
\begin{defn}\label{convexity radius}
{\em The convexity radius of a point $x\in M$ in a Riemannian
manifold $M$ is defined as the supremum in
$\overline{\mathbb{R}^{+}}$ of the numbers $r>0$ such that the
ball $B(x,r)$ is convex. We denote this supremum by $c(M,x)$. We
define the global convexity radius of $M$ as $c(M):=\inf\{ c(M,x)
: x\in M\}$.}
\end{defn}
\begin{rem}\label{continuity of the convexity radius}
{\em By Whitehead's theorem we know that $c(x,M)>0$ for every
$x\in M$. On the other hand, the function $x\mapsto c(x,M)$ is
continuous on $M$, see \cite[Corollary 1.9.10]{Klingenberg}.
Consequently, if $M$ is compact, then $c(M)>0$, that is, $M$ is
uniformly locally convex.}
\end{rem}

The notion of injectivity radius of a Riemannian manifold will
also play a role in the study of variational principles and
Hamilton-Jacobi equations. Let us recall its definition.
\begin{defn}
{\em We define the injectivity radius of a Riemannian manifold $M$
at a point $x\in M$ as the supremum in $\overline{\mathbb{R}^{+}}$
of the numbers $r>0$ such that $\exp_x$ is a $C^\infty$
diffeomorphism onto its image when restricted to the ball
$B(0_{x}, r)$. We denote this supremum by $i(M,x)$. The
injectivity radius of $M$ is defined by $i(M):=\inf\{ i(M,x) :
x\in M\}$.}
\end{defn}
\begin{rem}
{\em For a finite-dimensional manifold $M$, it can be seen that
$i(M,x)$ equals the supremum of the numbers $r>0$ such that
$\exp_x$ is injective when restricted to the ball $B(0_{x}, r)$,
see \cite{Klingenberg}. However, for infinite dimensional
manifolds it is not quite clear if this is always true.}
\end{rem}
\begin{rem}\label{continuity of the injectivity radius}
{\em By Theorem \ref{properties of exp} we know that $i(x,M)>0$
for every $x\in M$. On the other hand, it is well known that the
function $x\mapsto i(x,M)$ is continuous on $M$ \cite[Proposition
2.1.10]{Klingenberg}. Therefore, if $M$ is compact, then
$i(M)>0$.}
\end{rem}
We will also need to use the parallel translation of vectors along
geodesics. Recall that, for a given curve $\gamma: I\to M$, a
number  $t_{0}, t_{1}\in I$, and a vector $V_{0}\in
TM_{\gamma(t_{0})}$, there exists a unique parallel vector field
$V(t)$ along $\gamma(t)$ such that $V(t_{0})=V_{0}$. Moreover,
the mapping defined by $V_{0}\mapsto V(t)$ is a linear isometry
between the tangent spaces $TM_{\gamma(t_{0})}$ and
$TM_{\gamma(t)}$, for each $t\in I$. We denote this mapping by
$P_{t_{0}}^{t}=P_{t_{0}, \gamma}^{t}$, and we call it the parallel
translation from $TM_{\gamma(t_{0})}$ to $TM_{\gamma(t)}$ along
the curve $\gamma$.

The parallel translation will allow us to measure the length of
the ``difference" between vectors (or forms) which are in
different tangent spaces (or in duals of tangent spaces, that is,
fibers of the cotangent bundle), and do so in a natural way.
Indeed, let $\gamma$ be a minimizing geodesic connecting two
points $x, y\in M$, say $\gamma(t_{0})=x, \gamma(t_{1})=y$. Take
vectors $V\in TM_{x}$, $W\in TM_{y}$. Then we can define the
distance between $V$ and $W$ as the number
    $$
    \|W-P_{t_{0}, \gamma}^{t_{1}}(V)\|_{y}=
    \|V-P_{t_{1}, \gamma}^{t_{0}}(W)\|_{x}
    $$
(this equality holds because $P_{t_{0}}^{t_{1}}$ is a linear
isometry between the two tangent spaces, with inverse
$P_{t_{1}}^{t_{0}}$). Since the spaces $T^{*}M_{x}$ and $TM_{x}$
are isometrically identified by the formula $v=\langle v,
\cdot\rangle$, we can obviously use the same method to measure
distances between forms $\xi\in T^{*}M_{x}$ and $\eta\in
T^{*}M_{y}$ lying in different fibers of the cotangent bundle.

\medskip

Finally, let us consider some mean value theorems. The following
two results are easily deduced from the mean value theorem for
functions of one variable, but it will be convenient to state and
prove them for future reference.

\begin{thm}[Mean value theorem]\label{mean value}
Let $(M, g)$ be a Riemannian manifold, and $f:M\to\mathbb{R}$ a
Fr{\'e}chet differentiable mapping. Then, for every pair of points
$p, q\in M$ and every minimal geodesic path $\sigma:I\to M$
joining $p$ and $q$, there exists $t_{0}\in I$ such that
    $$
    f(p)-f(q)= d(p,q) df(\sigma(t_{0})) (\sigma'(t_{0}));
    $$
in particular
$|f(p)-f(q)|\leq\|df(\sigma(t_{0}))\|_{\sigma(t_{0})} d(p,q)$.
\end{thm}
\begin{proof}
Since $\sigma$ is a minimal geodesic we may assume $I=[0,
d(p,q)]$, $\|\sigma'(t)\|_{\sigma(t)}=1$ for all $t\in I$,
$\sigma(0)=q$, $\sigma(d(p,q))=p$. Consider the function
$h:I\to\mathbb{R}$ defined by $h(t)=f(\sigma(t))$. By applying the
mean value theorem to the function $h$ we get a point $t_{0}\in I$
such that
    $$
    f(p)-f(q)=h(d(p,q))-h(0)=h'(t_{0})(d(p,q)-0)=
    df(\sigma(t_{0}))(\sigma'(t_{0})) d(p,q),
    $$
and, since $\|\sigma'(t_{0})\|_{\sigma(t_{0})}=1$ and
$|df(\sigma(t_{0}))(\sigma'(t_{0}))|\leq
\|df(\sigma(t_{0}))\|_{\sigma(t_{0})}$, we also get that
$|f(p)-f(q)|\leq\|df(\sigma(t_{0}))\|_{\sigma(t_{0})} d(p,q)$.
\end{proof}
When the points cannot be joined by a minimal geodesic we have a
less accurate but quite useful result which tells us that every
function with a bounded derivative is Lipschitz with respect to
the g-distance on $M$. In fact this results holds even for
functions which take values in other Riemannian manifolds. For a
differentiable function between Riemannian manifolds $f:M\to N$,
we define the norm of the derivative $df(p)$ at a point $p\in M$
by
\begin{eqnarray*}
    & &\|df(p)\|_{p}:=\sup\{ \|df(p)(v)\|_{f(p)} : h\in TM_{p},
    \|v\|_{p}\leq 1\}=\\
    & &\sup\{ \zeta\big(df(p)(v)\big) : v\in TM_{p},
    \zeta\in T^{*}N_{f(p)}, \|v\|_{p}=1=\|\zeta\|_{f(p)}\}.
\end{eqnarray*}
\begin{thm}[Mean value inequality]\label{Mean value inequality}
Let $M, N$ be Riemannian manifolds, and $f:M\to N$ a Fr{\'e}chet
differentiable mapping. Assume that $f$ has a bounded derivative,
say $\|df(x)\|_{x}\leq C$ for every $x\in M$. Then $f$ is
C-Lipschitz, that is
    $$
    d_{N}\big(f(p), f(q)\big)\leq C d_{M}(p,q)
    $$
for all $p, q\in M$.
\end{thm}
\begin{proof}
Fix any two points $p, q\in M$. Take any $\varepsilon>0$. By
definition of $d(p,q)$, there exists a $C^1$ smooth path
$\gamma:[0, T]\to M$ with $\gamma(0)=q$, $\gamma(T)=p$, and
    $$
    L(\gamma)\leq d_{M}(p,q)+\frac{\varepsilon}{C};
    $$
as usual we may assume $\|\gamma'(t)\|_{\gamma(t)}=1$ for all
$t\in [0,T]=[0, L(\gamma)]$. By considering the path
$\beta(t):=f(\gamma(t))$, which joins the points $f(p)$ and $f(q)$
in $N$, and bearing in mind the definitions of $d_{N}(f(p),f(q))$
and the fact that $\|d\gamma(t)\|_{\gamma(t)}=1$ for all $t$, we
get
\begin{eqnarray*}
& & d_{N}\big(f(p),f(q)\big)\leq
L(\beta)=\int_{0}^{T}\|d\beta(t)\|_{\beta(t)}dt
=\int_{0}^{T}\|df(\gamma(t))(d\gamma(t))\|_{f(\gamma(t))}dt\leq\\
& &\int_{0}^{T}\|df(\gamma(t))\|_{\gamma(t)}dt\leq \int_{0}^{T}C
dt=C T\leq C\big(d_{M}(p,q)+\varepsilon/C\big)=C d_{M}(p,q)
+\varepsilon.
\end{eqnarray*}
We have shown that $d_{N}(f(p),f(q))\leq C d_{M}(p,q)+\varepsilon$
for every $\varepsilon>0$, which means that $d_{N}(f(p),f(q))\leq
C d_{M}(p,q)$.
\end{proof}

In Section 4 below we will generalize these mean value theorems
for the case of subdifferentiable or superdifferentiable functions
defined on Riemannian manifolds.

The preceding mean value theorem has a converse, which is
immediate in the case when $M$ and $N$ are Hilbert spaces, but
requires some justification in the setting of Riemannian
manifolds.
\begin{prop}\label{converse of the mean value inequality}
Let $M, N$ be Riemannian manifolds. If $f:M\to N$ is $K$-Lipschitz
(that is, $d_{N}(f(x),f(y))\leq K d_{M}(x,y)$ for all $x,y\in
M$), then $\|df(x)\|_{x}\leq K$ for every $x\in M$.
\end{prop}
\begin{proof}
Consider first the case when $N=\mathbb{R}$. Suppose that there
exists $x_{0}\in M$ with $\|df(x_{0})\|_{x_{0}}>K$. Take $v_{0}\in
TM_{x_{0}}$ so that $\|v_{0}\|_{x_{0}}=1$ and
$df(x_{0})(v_{0})>K$. Consider the geodesic
$\gamma(t)=\exp_{x_{0}}(tv_{0})$ defined for $|t|\leq r_{0}$ with
$r_{0}>0$ small enough. Define $F:[-r_{0},r_{0}]\to \mathbb{R}$
by $F(t)=f(\gamma(t))$. We have that $F'(0)=df(x_{0})(v_{0})>K$.
By the definition of $F'(0)$ we can find some $\delta_{0}\in
(0,r_{0})$ such that
    $$
    \frac{F(t)-F(0)}{t}>K \textrm{ if } |t|\leq\delta_{0}.
    $$
Taking $t_{1}=-\delta_{0}$, $t_{2}=\delta_{0}$ we get
$F(t_{1})-F(0)<K t_{1}$ and $F(t_{2})-F(0)>K t_{2}$, hence, by
summing, $$F(t_{2})-F(t_{1})>K(t_{2}-t_{1}).$$ If we set
$x_{1}=\gamma(t_{1})$, $x_{2}=\gamma(t_{2})$ this means that
    $$
    f(x_{2})-f(x_{1})>K(t_{2}-t_{1})=K d(x_{2}, x_{1}),
    $$
which contradicts the fact that $f$ is $K$-Lipschitz.

Now let us consider the general case when the target space is a
Riemannian manifold $N$. Suppose that $\|df(x_{0})\|_{x_{0}}>K$
for some $x_{0}\in M$. Then there are $\zeta_{0}\in
T^{*}N_{f(x_{0})}$ and $v_{0}\in TM_{x_{0}}$ with
$\|v_{0}\|_{x_{0}}=1=\|\zeta_{0}\|_{f(x_{0})}$ and such that
$K<\|df(x_{0})\|_{x_{0}}=\zeta_{0}\big(df(x_{0})(v_{0})\big)$.
Take $s_{0}>0$ and $\varepsilon>0$ small enough so that
$\exp_{f(x_{0})}^{-1}:B(f(x_{0}),s_{0})\to B(0_{f(x_{0})},
s_{0})$ is a $(1+\varepsilon)$-Lipschitz diffeomorphism and
$K<(1+\varepsilon)K<\|df(x_{0})\|_{x_{0}}$. Now take $r_{0}>0$
small enough so that $f(B(x_{0}, r_{0}))\subset B(f(x_{0}),
s_{0})$, and define the composition
    $$
    g:B(x_{0},r_{0})\to\mathbb{R}, \,\,
    g(x)=\zeta_{0}\big(\exp_{f(x_{0})}^{-1}(f(x))\big).
    $$
It is clear that $g$ is $(1+\varepsilon) K$-Lipschitz. But, since
$d\exp_{f(x_{0})}^{-1}(f(x_{0}))$ is the identity, we have that
\begin{eqnarray*}
& &dg(x_{0})(v_{0})=\zeta_{0}\big(d\exp_{f(x_{0})}^{-1}(f(x_{0}))
    (df(x_{0})(v_{0}))\big)=\\
&
&\zeta_{0}\big(df(x_{0})(v_{0})\big)=\|df(x_{0})\|_{x_{0}}>(1+\varepsilon)K,
\end{eqnarray*}
and this contradicts the result we have just proved for the case
$N=\mathbb{R}$.
\end{proof}
\medskip


\section{Almost-critical-point-spotting results}

As said in the introduction, in infinite dimensions one cannot
generally hope to find any critical point for a given smooth
function, whatever its shape, so one has to make do with almost
critical points.

\medskip

\noindent {\bf An approximate Rolle's theorem}

We begin with an approximate version of Rolle's theorem which
holds in every Riemannian manifold (even though it is
infinite-dimensional) and ensures that every differentiable
function which has a small oscillation on the boundary of an open
set whose closure is complete has an almost critical point.

\begin{thm}[Approximate Rolle's
theorem]\label{Rolle for manifolds} Let $(M,g)$ be a Riemannian
manifold, $U$ an open subset of $M$ such that $\overline{U}$ is
complete and bounded with respect to the $g$-distance, and
$p_{0}\in M$, $R>0$ be such that $B_{g}(p_{0},R)\subseteq
\overline{U}$. Let $f:\overline{U}\longrightarrow\mathbb{R}$ be a
continuous function which is differentiable on $U$. Then:
\begin{enumerate}
\item If $\sup f(U)>\sup f(\partial U)$ then, for every $r>0$
there exists $q\in U$ such that $\|df(q)\|_{q}\leq r$.
\item If $\inf f(U)<\inf f(\partial U)$ then, for every $r>0$
there exists $q\in U$ such that $\|df(q)\|_{q}\leq r$.
\item If $f(\overline{U})\subseteq[-\varepsilon,\varepsilon]$ for
some $\varepsilon>0$, then there exists $q\in U$ such that
$\|df(q)\|_{q}\leq\varepsilon/R$.
\end{enumerate}
\end{thm}
\begin{cor}
Let $(M,g)$ be a complete Riemannian manifold, $U$ a bounded open
subset of $M$, and $p_{0}\in M$, $R>0$ be such that
$B_{g}(p_{0},R)\subseteq \overline{U}$, $\varepsilon>0$. Suppose
that $f(\partial U)\subseteq[-\varepsilon,\varepsilon]$. Then
there exists some $q\in U$ such that
$\|df(q)\|_{q}\leq\varepsilon/R$.
\end{cor}

\medskip

To prove Theorem \ref{Rolle for manifolds} we begin with a simple
lemma.
\begin{lem}\label{decreasing principle}
Let $(M, g)$ be a Riemannian manifold, and
$f:M\longrightarrow\mathbb{R}$ be a differentiable function on
$M$. Suppose that $\|df(p)\|_{p}>\varepsilon>0$. Then there exist
a number $\delta>0$ and two $C^1$ paths
$\alpha,\beta:[0,\delta]\to M$, parameterized by arc length, such
that
    $$
    f(\alpha(t))<f(p)-\varepsilon t, \, \text{ and }\,
    f(\beta(t))>f(p)+\varepsilon t,
    $$
for all $t\in (0,\delta]$.
\end{lem}
\begin{proof}
Let us show the existence of such a path $\alpha$ (a required path
$\beta$ can be obtained in a similar manner). Since
$\|df(p)\|_{p}>\varepsilon$, there exists $h\in TM_{p}$ so that
$\|h\|_{p}=1$ and $df(p)(h)<-\varepsilon$. Then (by the
characterization of the tangent space $TM_{p}$ as the set of
derivatives of all smooth paths passing through $p$) we can choose
a $C^1$ path $\alpha:[0,r]\to M$, parameterized by arc length,
such that $$\frac{d\alpha}{dt}(0)=h, \text{ and } \alpha(0)=p.$$
Define the function $F:[0,r]\to\mathbb{R}$ by
    $F(t)=f(\alpha(t)).$ We have that
    $$
    F'(s)=df(\alpha(s))(\frac{d\alpha}{dt}(s))
    $$
for all $s\in [0,r]$. In particular, for $s=0$, we have that
    $F'(0)=df(p)(h)<\varepsilon,$ and therefore there exists some
$\delta>0$ such that
    $$
    \frac{F(t)-F(0)}{t}<-\varepsilon
    $$
for all $t\in (0,\delta]$. This means that
    $f(\alpha(t))<f(p)-\varepsilon t$ for all $t\in (0,\delta]$.
\end{proof}

We will also make use of the following version of Ekeland's
Variational Principle (see \cite{Ekeland1} for a proof).
\begin{thm}[Ekeland's Variational Principle]\label{Ekeland}
Let $X$ be a complete metric space, and let $f:X\longrightarrow
[-\infty, \infty)$ be a proper upper semicontinuous function which
is bounded above. Let $\varepsilon
>0$ and $x_{0}\in X$ such that $f(x_{0})>\sup\{f(x) : x\in X\}
-\varepsilon$. Then for every $\lambda$ with $0<\lambda<1$ there
exists a point $z\in Dom(f)$ such that:
\begin{enumerate}
\item[(i)] $\lambda d(z,x_{0})\leq f(z)-f(x_{0})$
\item[(ii)] $d(z, x_{0}) <\varepsilon/\lambda$
\item[(iii)] $\lambda d(x,z)+f(z)> f(x)$ whenever $x\neq z$.
\end{enumerate}
\end{thm}

\medskip

\begin{center}
{\bf Proof of Theorem \ref{Rolle for manifolds}.}
\end{center}

\noindent{\bf Case 1}: Let $\eta =\sup f(U) -\sup f(\partial U
)>0$. Define $X=(\overline{U}, d_{g})$, which is a complete metric
space. Let $n>1$ be large enough so that $\overline{U}\subset
B_{g}(p_{0}, n)$, and set $\lambda=\min\{\eta/8n, r\}>0$. Observe
that, since the diameter of $U$ is less than or equal to $2n$, we
have that $\lambda d(x,y)\leq\eta/4$ for all $x,y\in\overline{U}$.
Now, according to Ekeland's Variational Principle \ref{Ekeland},
there exists $q\in \overline{U}$ such that
     $$
     f(y)\leq f(q)+\lambda d(y,q) \text{ for
all } y\in X. \eqno(1)
      $$ In fact, it must be $q\in U$: if $q\in\partial U$
then, taking $a$ such that $f(a)\geq \sup f(U)-\eta/4$ we would
get $$\sup f(U) - \eta/2=(\sup f(U)-\eta/4)-\eta/4 \leq f(a)
-\lambda d(a,q) \leq f(q) \leq \sup f(\partial U),$$ a
contradiction.

We claim that $\|df(q)\|_{q}\leq \lambda\leq r$. Indeed, assume
that $\|df(q)\|_{q}>\lambda$. Then, according to Lemma
\ref{decreasing principle}, there would exist a $C^1$ path
$\beta$, parameterized by arc length, such that $\beta(0)=q$ and
    $$
    f(\beta(t))>f(q)+\lambda t \eqno(2)
    $$
for all $t>0$ small enough. By combining $(1)$ and $(2)$, we would
get that
    $$
   f(q)+\lambda t< f(\beta(t))\leq f(q)+\lambda d(\beta(t),q)
   \leq f(q)+\lambda L(\beta_{|_{[0,t]}})=f(q)+\lambda t
   $$
if $t>0$ is small enough; but this is a contradiction.

\noindent{\bf Case 2}: consider the function $-f$ and apply Case
(1).

\noindent{\bf Case 3}: We will consider two situations.

\noindent{\bf Case 3.1}: Suppose that $f(p_{0})\neq 0$. We may
assume that $f(p_{0})<0$ (the case $f(p_{0})>0$ is analogous).
Define $\lambda=\varepsilon/R$. According to Ekeland's Variational
Principle, there exists $q\in\overline{U}$ such that
\begin{itemize}
\item[{(i)}] $d(p_{0}, q)\leq\frac{1}{\lambda}(f(p_{0})-f(q))\leq
\frac{1}{\lambda}(f(p_{0})+\varepsilon)<R$, and
\item[{(ii)}] $f(q)<f(y)+\lambda d(y,q)$ if $y\neq q$.
\end{itemize}
The first property tells us that
$q\in\text{int}B_{g}(p_{0},R)\subseteq U$. And, by using Lemma
\ref{decreasing principle} as in Case 1, it is immediately seen
that the second property implies that
$\|df(q)\|_{q}\leq\lambda=\varepsilon/R$.

\noindent{\bf Case 3.2}: Suppose finally that $f(p_{0})=0$. We
may assume that $\|df(p_{0})\|_{p_{0}}>\varepsilon/R$ (otherwise
we are done). By Lemma \ref{decreasing principle}, there exist
$\delta>0$ and a $C^1$ path $\alpha$ in $U$ such that
$$f(\alpha(t))<f(p_{0})-\frac{\varepsilon}{R} t$$ if
$0<t\leq\delta$. Define $x_{0}=\alpha(\delta)\in
B_{g}(p_{0},\delta)$. We have that
    $$
    f(x_{0})<f(p_{0})-\frac{\varepsilon}{R}\delta=-\frac{\varepsilon}{R}\delta<0.
    $$
By applying again Ekeland's Variational Principle with
$\lambda=\varepsilon/R$ we get a point $q\in\overline{U}$ such
that
\begin{itemize}
\item[{(i)}] $d(q,x_{0})\leq\frac{f(x_{0})+\varepsilon}{\varepsilon}<
\frac{-\varepsilon\delta/R+\varepsilon}{\varepsilon/R}=R-\delta$,
and
\item[{(ii)}] $f(q)<f(y)+\frac{\varepsilon}{R} d(y, q)$ for all $y\neq
q$.
\end{itemize}
Now, $(i)$ implies that $d(q, p_{0})\leq d(q, x_{0})+d(x_{0},
p_{0})<R-\delta+\delta=R$, that is, $q\in
\text{int}B_{g}(p_{0},R)\subseteq U$. And, as above, bearing in
mind Lemma \ref{decreasing principle}, $(ii)$ implies that
$\|df(q)\|_{q}\leq\varepsilon/R$. \qed
\begin{rem}\label{completeness is necessary}
{\em If $\overline{U}$ is not complete the result is obviously
false: consider for instance $M=(-1,1)\subset\mathbb{R}$,
$U=(0,1)$, $\partial U=\{0\}$, $f(x)=x$. On the other hand, the
estimate $\varepsilon/R$ is sharp, as this example shows:
$M=\mathbb{R}$, $U=(-1,1)$, $f(x)=x$, $R=1$, $p_{0}=0$,
$\varepsilon=1$.}
\end{rem}

\medskip

\noindent {\bf A smooth variational principle}

Now we turn our attention to perturbed minimization principles on
Riemannian manifolds. Of course, since every Riemannian manifold
is a metric space, Ekeland's variational principle quoted above
holds true and is very useful in this setting: every lower
semicontinuous function can be perturbed with a function whose
shape is that of an almost flat cone in such a way that the
difference attains a global minimum. But sometimes, especially
when one wants to build a good theory of subdifferentiability, one
needs results ensuring that the perturbation of the function is
smooth, that is, one needs to replace that cone with a smooth
function which is arbitrarily small and has an arbitrarily small
Lipschitz constant. This is just what the Deville-Godefroy-Zizler
Smooth Variational Principle does in those Banach spaces having
$C^1$ smooth Lipschitz bump functions; see \cite{DGZ}.

Unfortunately, the main ideas behind the proof of this Variational
Principle in the case of Banach spaces cannot be transferred to
the setting of Riemannian manifolds in full generality. One has to
impose some restriction on the structure of the manifold in order
that those ideas work. That is why we need the following
definition.

\begin{defn}\label{uniformly bumpable manifold}
We will say that a Riemannian manifold $M$ is {\em uniformly
bumpable} provided there exist numbers $R>1$ (possibly large) and
$r>0$ (small) such that for every $p\in M$, $\delta\in (0, r)$
there exists a $C^1$ smooth function $b:M\to\mathbb [0,1]$ such
that:
\begin{enumerate}
\item $b(p)=1$
\item $b(x)=0$ if $d(x,p)\geq\delta$
\item $\sup_{x\in M}\|db(x)\|_{x}\leq R/\delta$.
\end{enumerate}
\end{defn}
\begin{rem}
{\em It is easy to see that every Riemannian manifold $M$ is {\em
bumpable}, in the sense that for every $p\in M$, $\delta>0$, there
exists a smooth bump function $b:M\to [0,1]$ with $b(p)=1$,
$b(x)=0$ for $x\notin B(p,\delta)$, and $b$ is Lipschitz, that is
$\sup_{x\in M}\|db(x)\|_{x}<\infty$. However it is not quite clear
which Riemannian manifolds are uniformly bumpable. Of course every
Hilbert space is uniformly bumpable, and there are many other
natural examples of uniformly bumpable Riemannian manifolds. In
fact we do not know of any Riemannian manifold which is not
uniformly bumpable.}
\end{rem}
\begin{op}
{\em Is every Riemannian manifold uniformly bumpable? If not,
provide useful characterizations of those Riemannian manifolds
which are uniformly bumpable.}
\end{op}
The following Proposition provides some sufficient conditions for
a Riemannian manifold to be uniformly bumpable:  it is enough that
$\exp_{x}$ is a diffeomorphism and preserves radial distances
when restricted to balls of a fixed radius $r>0$. This is always
true when $M$ is uniformly locally convex and has a strictly
positive injectivity radius.
\begin{prop}\label{sufficient conditions for M to be
uniformly bumpable} Let $M$ be a Riemannian manifold. Consider the
following six conditions:
\begin{enumerate}
\item $M$ is compact.
\item $M$ is finite-dimensional, complete, and has a strictly positive
injectivity radius $i(M)$.
\item $M$ is uniformly locally convex and has a strictly positive injectivity radius.
\item There is a constant $r>0$ such that for every
$x\in M$ the mapping $\exp_{x}$ is defined on $B(0_{x}, r)\subset
TM_{x}$ and provides a $C^\infty$ diffeomorphism
    $$
    \exp_{x}:B(0_{x}, r)\to B(x,r),
    $$
and the distance function is given here by the expression
    $$
    d(y,x)=\|\exp_{x}^{-1}(y)\|_{x} \, \textrm{ for all }\, y\in
    B(x,r).
    $$
\item There is a constant $r>0$ such that for every
$x\in M$ the distance function to $x$, $y\in M\mapsto d(y,x)$, is
$C^\infty$ smooth on the punctured ball $B(x,r)\setminus\{x\}$.
\item $M$ is uniformly bumpable.
\end{enumerate}
Then $(1)\implies (2)\implies (3)\iff (4)\implies (5)\implies
(6)$.
\end{prop}
\begin{proof}
\noindent $(1)\implies (2)$ is a trivial consequence of Remark
\ref{continuity of the injectivity radius}.

\medskip

\noindent $(2)\implies (3)$: In \cite[Chapter 2]{Klingenberg}, the
injectivity radius of a point $x\in M$ is characterized as the
distance from $x$ to the cut locus $C(x)$ of $x$. Hence, for every
$r>0$ with $r< i(M)$ and every $x\in M$ it is clear that
$\exp_{x}$ is a diffeomorphism and preserves radial distances when
restricted to balls of a fixed radius $r>0$ in the tangent space
$TM_{x}$, and $M$ is uniformly locally convex. See Theorems 2.1.14
and 2.1.12 of \cite{Klingenberg}.

\medskip

\noindent $(3)\implies (4)$: Since $i(M)>0$, we know that there is
some $r_{1}>0$ such that $\exp_{x}$ is a diffeomorphism onto its
image when restricted to the ball $B(0_{x}, r_{1})$, for all $x\in
M$. The fact that $M$ is uniformly locally convex clearly implies
that there is some $r_{2}>0$ such that
    $$
    d(y,x)=\|\exp_{x}^{-1}(y)\|_{x} \, \textrm{ for all }\, y\in
    B(x,r_{2}).
    $$
We may obviously assume that $r_{1}=r_{2}:=r$. In particular
$\exp_{x}$ maps $B(0_{x},r)$ onto $B(x,r)$.

\medskip

\noindent $(4)\implies (3)$ is obvious.

\medskip

\noindent $(4)\implies (5)$ is trivial, since $\exp_{x}^{-1}$ is a
$C^\infty$ diffeomorphism between those balls, $\|\cdot\|_{x}$ is
$C^\infty$ smooth on $TM_{x}\setminus\{0_{x}\}$, and
$d(y,x)=\|\exp_{x}^{-1}(y)\|_{x}$ for all $y\in B(x,r)$.

\medskip

\noindent $(5)\implies (6)$: Assume that the distance function
$y\mapsto d(y,x)$ is $C^\infty$ smooth on $B(x,r)\setminus\{x\}$.
Let $\theta:\mathbb{R}\to [0,1]$ be a $C^{\infty}$ smooth
Lipschitz function such that $\theta^{-1}(1)=(-\infty, 1/3]$ and
$\theta^{-1}(0)=[1, \infty)$. For a given point $x\in M$ and a
number $\delta\in (0, r)$, define $b:M\to [0,1]$ by
    $$
    b(y)=\theta\big(\frac{1}{\delta}d(y,x)\big).
    $$
Taking into account the fact that the distance function $y\mapsto
d(y,x)$ is 1-Lipschitz and therefore the norm of its derivative
is everywhere bounded by $1$ (see Proposition \ref{converse of the
mean value inequality}), it is easy to check that $b$ satisfies
conditions 1-2-3 of Definition \ref{uniformly bumpable manifold},
for a constant $R=\|\theta'\|_{\infty}>1$ that only depends on the
real function $\theta$, but not on the point $x\in M$.
\end{proof}
\begin{rem}
{\em The condition that $M$ has a strictly positive injectivity
radius is not necessary in order that $M$ is uniformly bumpable,
as the following example shows. Let $M$ be the surface of
$\mathbb{R}^{3}$ defined by the equation $z=1/(x^{2}+y^{2})$,
$(x,y)\neq (0,0)$, with the natural Riemannian structure inherited
from $\mathbb{R}^{3}$. Then $i(M)=0$, but, as is not difficult to
see, $M$ is uniformly bumpable.}
\end{rem}
The following theorem is the natural extension of the
Deville-Godefroy-Ziz\-ler Smooth Variational Principle to
Riemannian manifolds which are uniformly bumpable. Recall that a
function $F:M\to\mathbb{R}\cup\{+\infty\}$ is said to attain a
strong minimum at $p$ provided $F(p)=\inf_{x\in M}F(x)$ and
$\lim_{n\to\infty}d(p_{n},p)=0$ whenever $(p_{n})$ is a minimizing
sequence (that is, if \, $\lim_{n\to\infty}F(p_{n})=F(p)$).
\begin{thm}[DGZ Smooth Variational Principle]\label{Smooth Variational Principle}
Let $(M,g)$ be a complete Riemannian manifold which is uniformly
bumpable, and let $F:M\longrightarrow (-\infty,+\infty]$ be a
lower semicontinuous function that is bounded below,
$F\not\equiv+\infty$. Then, for every $\delta>0$ there exists a
bounded $C^1$ smooth function $\varphi:M\longrightarrow\mathbb{R}$
such that:
\begin{enumerate}
\item $F-\varphi$ attains its strong minimum in $M$,
\item $\|\varphi\|_{\infty}:=\sup_{p\in M}|\varphi(p)|<\delta$, and
$\|d\varphi\|_{\infty}:=\sup_{p\in M}\|d\varphi(p)\|_{p}<\delta$.
\end{enumerate}
\end{thm}
\begin{rem}
{\em The assumption that $M$ is complete is necessary here, as the
following trivial example shows: $M=(-1,1)\subset\mathbb{R}$,
$f(x)=x$.}
\end{rem}
We will split the proof of Theorem \ref{Smooth Variational
Principle} into three lemmas. In the sequel $B(x,r)$ denotes the
open ball of center $x$ and radius $r$ in the metric space $M$,
and $B(\varphi, r)$ stands for the open ball of center $\varphi$
and radius $r$ in the Banach space $Y$.
\begin{lem}\label{the superlemma}
Let $M$ be a complete metric space, and $(Y,\|\cdot\|)$ be a
Banach space of real-valued bounded and continuous functions on
$M$ satisfying the following conditions:
\begin{enumerate}
\item $\|\varphi\|\geq\|\varphi\|_{\infty}=\sup\{|\varphi(x)|:x\in
M\}$ for every $\varphi\in Y$.
\item There are numbers $C>1, r>0$ such that for every $p\in M$,
$\varepsilon>0$ and $\delta\in (0, r)$ there exists a function
$b\in Y$ such that $b(p)=\varepsilon$, $\|b\|_{Y}\leq
C\varepsilon (1+1 /\delta)$, and $b(x)=0$ if $x\not\in
{B(p,\delta)}$.
\end{enumerate}
Let $f:M\rightarrow \mathbb{R}\cup\{+\infty\}$ be a lower
semicontinuous function which is bounded below and such that
$Dom(f)=\{x\in M | f(x)<+\infty\}\neq\emptyset$. Then, the set $G$
of all the functions $\varphi\in Y$ such that $f+\varphi$ attains
a strong minimum in $M$ contains a $G_{\delta}$ dense subset of
$Y$.
\end{lem}
\begin{proof}
Take a number $N\in\mathbb{N}$ such that $N\geq 1/r$, and for
every $n\in\mathbb{N}$ with $n\geq N$, consider the set
$$U_n=\{\varphi\in Y | \, \exists \, x_0 \in M \,
:\, (f+\varphi)(x_0)<\inf\{(f+\varphi)(x)|x\in M\backslash
B(x_{0},\frac{1}{n})\}\}.$$ Let us see that $U_n$ is an open dense
subset of $Y$. Indeed,

\noindent $\bullet$ $U_n$ is open. Take $\varphi\in U_n$. By the
definition of $U_n$ there exists $x_{0}\in M$ such that
$(f+\varphi)(x_0)<\inf\{(f+\varphi)(x)|x\in M\backslash
B(x_{0},\frac{1}{n})\}$. Set $2\rho=\inf\{(f+\varphi)(x)|x\in
M\backslash B(x_{0},\frac{1}{n})\}-(f+\varphi)(x_0)>0$. Then,
since $\|\cdot\|_Y\geq \|\cdot\|_{\infty}$, we get that
$B_{Y}(\varphi, \rho)\subset B_{\infty}(\varphi,\rho)\subset U_n$.

\noindent $\bullet$ $U_n$ is dense in $Y$. Take $\varphi\in Y$ and
$\varepsilon>0$. Since $f+\varphi$ is bounded below there exists
$x_0 \in M$ such that $(f+\varphi)(x_0)<
\inf\{(f+\varphi)(x)|x\in M\}+\varepsilon$. Set now
$\delta=1/n<r$, and use condition (2) to find a function $b\in Y$
such that $b(x_0)=\varepsilon$, $\|b\|_{Y}\leq
C(n+1)\varepsilon$, and $b(x)=0$ for $x\not\in
B(x_{0},\frac{1}{n})$. Then
$(f+\varphi)(x_0)-b(x_0)<\inf\{(f+\varphi)(x)|x\in M\}$ and, if we
define $h=-b$, we have
    $$(f+\varphi+h)(x_0)<\inf\{(f+\varphi)(x)|x\in M\}\leq
    \inf\{(f+\varphi)(x)|x\not\in
B(x_0,\frac{1}{n})\}.$$ Since $\inf\{(f+\varphi)(x)|x\not\in
B(x_0,\frac{1}{n})\}=\inf\{(f+\varphi+h)(x)|x\not\in
B(x_0,\frac{1}{n})\}$, it is obvious that the above inequality
implies that $\varphi+h\in U_n$. On the other hand, we have
$\|h\|_{Y}\leq C(n+1)\varepsilon$. Since $C$ and $n$ are fixed and
$\varepsilon$ can be taken to be arbitrarily small, this shows
that $\varphi\in\overline{U_{n}}$, and $U_{n}$ is dense in $Y$.

Therefore we can apply Baire's theorem to conclude that the set
$G=\bigcap_{n=N}^{\infty} U_n$ is a $G_\delta$ dense subset of
$Y$. Now we must show that if $\varphi\in G$ then $f+\varphi$
attains a strong minimum in $M$. For each $n\geq N$, take $x_n \in
M$ such that $(f+\varphi)(x_n)<\inf\{(f+\varphi)(x)|x\not\in
B(x_n,\frac{1}{n})\}$. Clearly, $x_k\in B(x_n,\frac{1}{n})$ if
$k\geq n$, which implies that $(x_n)_{n=N}^{\infty}$ is a Cauchy
sequence in $M$ and therefore converges to some $x_0 \in M$.
Since $f$ is lower semicontinuous and
$\bigcap_{n=N}^{\infty}B(x_{0}, 1/n)=\{x_{0}\}$, we get
\begin{eqnarray*}
&&(f+\varphi)(x_{0})\leq\liminf(f+\varphi)(x_n)\leq
\liminf[\inf\{(f+\varphi)(x)|x\in M\setminus
B(x_0,\frac{1}{n})\}]\\ &&=\inf\{\inf\{(f+\varphi)(x)|x\in
M\setminus B(x_{0},\frac{1}{n})\}
: n\in\mathbb{N}, n\geq N\}\\
&&=\inf\{(f+\varphi)(x)|x\in M\setminus \{x_0\}\},
\end{eqnarray*}
which means that $f+\varphi$ attains a global minimum at $x_0\in
M$.

Finally, let us check that in fact $f+\varphi$ attains a strong
minimum at the point $x_{0}$. Suppose $\{y_n\}$ is a sequence in
$M$ such that $(f+g)(y_n)\rightarrow (f+g)(x_0)$ and $(y_n)$ does
not converge to $x_0$. We may assume $d(y_n, x_0)\geq \varepsilon$
for all $n$. Bearing in mind this inequality and the fact that
$x_{0}=\lim x_{n}$, we can take $k\in\mathbb{N}$ such that
$d(x_k, y_n)>\frac{1}{k}$ for all $n$, and therefore
    $$(f+\varphi)(x_0)\leq
(f+\varphi)(x_k)<\inf\{(f+\varphi)(x) | x\notin B(x_{k},
\frac{1}{k}) \}\leq(f+\varphi)(y_n)$$ for all $n$, which
contradicts the fact that $(f+\varphi)(y_n)\rightarrow
(f+\varphi)(x_0)$.
\end{proof}
\begin{lem}
Let $M$ be a uniformly bumpable Riemannian manifold. Then there
are numbers $C>1, r>0$ such that for every $p\in M$,
$\varepsilon>0$ and $\delta\in (0, r)$ there exists a $C^1$ smooth
function $b:M\to [0,\varepsilon]$ such that:
\begin{enumerate}
\item $b(p)=\varepsilon= \|b\|_{\infty}:=\sup_{x\in M}|b(x)|$.
\item $\|db\|_{\infty}:=\sup_{x\in M}\|db(x)\|_{x}\leq
C\varepsilon/\delta$.
\item $b(x)=0$ if $x\not\in {B(p,\delta)}$.
\end{enumerate}
In particular, $\max\{\|b\|_{\infty}, \|db\|_{\infty}\}\leq
C\varepsilon (1+1/\delta)$.
\end{lem}
\begin{proof}
The definition of uniformly bumpable manifold provides such $b$ in
the case when $\varepsilon=1$. If $\varepsilon\neq 1$, it is
enough to consider $b_{\varepsilon}=\varepsilon b$.
\end{proof}
\begin{lem}
Let $(M,g)$ be a complete Riemannian manifold. Then the vector
space $Y=\{\varphi:M\rightarrow\mathbb{R} \, | \, \varphi \textrm{
{\em is $C^1$ smooth, bounded and Lipschitz}}\}$, endowed with
the norm
$\|\varphi\|_{Y}=\max\{\|\varphi\|_{\infty},\|d\varphi\|_{\infty}
\}$, is a Banach space.
\end{lem}
\begin{proof}
It is obvious that $(Y, \|\cdot\|_{Y})$ is a normed space. We
only have to show that $Y$ is complete. Let $(\varphi_n)$ be a
Cauchy sequence with respect to the norm $\|\cdot\|_Y$. Since the
uniform limit of a sequence of continuous mappings between metric
spaces is continuous, it is obvious that $(\varphi_n)$ uniformly
converges to a continuous function $\varphi:M\to\mathbb{R}$.
Since $T^{*}M_{x}$ is a complete normed space for each $x\in M$,
it is also clear that $(d\varphi_n)$ converges to a function
$\psi:M\to T^{*}M$ defined by
    $$
    \psi(x)=\lim_{n\to\infty}d\varphi_{n}(x)
    $$
(where the limit is taken in $TM_{x}$ for each $x\in M$). Let us
see that $\psi=d\varphi$. Take $p\in M$. From Theorem
\ref{properties of exp} we know that there exists some $r>0$
(depending on $p$) such that the exponential mapping is defined on
$B(0_{p}, r)\subset TM_{p}$ and gives a diffeomorphism
$\exp_{p}:B(0_{p}, r)\to B(p,r)$ such that the derivatives of
$\exp_{p}$ and its inverse $(\exp_{p})^{-1}$ are bounded by $2$ on
$B(0_{p}, r)$ and $B(p,r)$ respectively; in particular $\exp_{p}$
provides a bi-Lipschitz diffeomorphism between these balls. We
denote $\widetilde{\varphi}(h)=(\varphi\circ\exp_p)(h)$, for $h\in
B(0_{p}, r)$, and $\widetilde{\varphi}(0_{p})=d\varphi(p)\circ
(d\exp_{p}(0_{p}))=d\varphi(p)$. We have
\begin{eqnarray*}
& &\big|\frac{\widetilde{\varphi}(h)-\widetilde{\varphi}(0)-
\psi(p)(h)}{\|h\|}\big|=
\big|\frac{\widetilde{\varphi}(h)-\widetilde{\varphi}(0)}{\|h\|}-
\psi(p)(\frac{h }{\|h\|})\big|\leq \\ &
&\big|\frac{\widetilde{\varphi}(h)-\widetilde{\varphi}(0)-(\widetilde{\varphi}_n
(h)-\widetilde{\varphi}_n
(0))}{\|h\|}\big|+\big|\frac{\widetilde{\varphi}_n
(h)-\widetilde{\varphi}_n (0)}{\|h\|}-d\widetilde{\varphi}_n
(0)(\frac{h}{\|h\|})\big|+ \\ & &\big|(d\widetilde{\varphi}_n
(0)-\psi(p))(\frac{h}{\|h\|})\big|. \hspace{7.8cm} (1)
\end{eqnarray*}
Let us first consider the expression
$|\frac{\widetilde{\varphi}(h)-\widetilde{\varphi}(0)-(\widetilde{\varphi}_n
(h)-\widetilde{\varphi}_n (0))}{\|h\|}|$. By applying the mean
value inequality theorem we get
\begin{eqnarray*}
& &|\widetilde{\varphi}_m(h)-\widetilde{\varphi}_m(0)-
(\widetilde{\varphi}_n (h)-\widetilde{\varphi}_n
(0))|\leq\sup_{x\in B(0_{p},r)} \|d\widetilde{\varphi}_{m}(x)
-d\widetilde{\varphi}_{n}(x)\|_{p}\|h\|_{p}\leq\\ &
&2\|d\varphi_{m}-d\varphi_{n}\|_{\infty}\|h\|_{p}.
\end{eqnarray*}
Since $(\varphi_{n})$ is a Cauchy sequence in $Y$ we deduce that
for every $\varepsilon>0$ there exists $n_0\in\mathbb{N}$ such
that $|\widetilde{\varphi}_m (h)-\widetilde{\varphi}_m
(0)-(\widetilde{\varphi}_n (h)-\widetilde{\varphi}_n
(0))|<(\varepsilon/3)\|h\|$ whenever $m,n\geq n_0$ so, by letting
$m\rightarrow\infty$ we get that
$|\widetilde{\varphi}(h)-\widetilde{\varphi}(0)-(\widetilde{\varphi}_n
(h)-\widetilde{\varphi}_n (0))|<(\varepsilon/3)\|h\|$ if $n\geq
n_0$.

On the other hand, the term $|(d\widetilde{\varphi}_n
(0)-\widetilde{\psi}(p))(\frac{h}{\|h\|})|$ in the right side of
inequality $(1)$ above is less than $\varepsilon/3$ when $n$ is
large enough; we may assume this happens if $n\geq n_{0}$.

Finally, if we fix $n=n_{0}$, the term
$|\frac{\widetilde{\varphi}_{n_{0}} (h)-\widetilde{\varphi}_n
(0)}{\|h\|}-d\widetilde{\varphi}_{n_{0}} (0)(\frac{h}{\|h\|})|$
can be made to be less than $\varepsilon/3$ if $\|h\|$ is small
enough, say $\|h\|\leq\delta$.

By combining these estimations we get that, for $n=n_{0}$, the
left side of inequality $(1)$ is less than $\varepsilon$ if
$\|h\|\leq\delta$. This shows that $\widetilde{\varphi}$ is
differentiable at $p$, with $d\widetilde{\varphi}(0_{p})=\psi(p)$.
Hence $\varphi$ is differentiable at $p$, with
$d\varphi(p)=\psi(p)$.

To conclude that $Y$ is a Banach space it only remains to check
that $d\varphi=\psi$ is continuous and bounded. Take
$\varepsilon>0$. Since $(\varphi_{n})$ is a Cauchy sequence in
$Y$, there exists $n_{0}\in\mathbb{N}$ such that
$\|d\varphi_{n}(y)-d\varphi_{m}(y)\|_{y}\leq\varepsilon$ for all
$y\in M$ provided $n,m\geq n_{0}$. By letting $m\to\infty$ we
deduce that $\|d\varphi_{n}(y)-\psi(y)\|_{y}\leq\varepsilon$ for
all $y\in M$, if $n\geq n_{0}$. That is, we have
    $$
    \lim_{n\to\infty}\|d\varphi_{n}-d\varphi\|_{\infty}=0.
    $$
In particular, this implies that $\|d\varphi\|_{\infty}<\infty$,
that is, $\varphi$ is Lipschitz. Now we can show $\psi=d\varphi$
is continuous. Take any $p\in M$. As above, there exists $r>0$
such that $\exp_{p}:B(0_{p}, r)\to B(p,r)$ is a 2-Lipschitz
diffeomorphism, and so is the inverse $\exp_{p}^{-1}$. Define
$\widetilde{\varphi}=\varphi\circ\exp_{p}:B(0_{p},r)\to\mathbb{R}$.
In order to see that $d\varphi$ is continuous at $p$ it is enough
to see that  $d\widetilde{\varphi}$ is continuous at $0_{p}$. By
applying the mean value inequality we have that
\begin{eqnarray*}
& &\|d\widetilde{\varphi}(x)-d\widetilde{\varphi}(0)\|_{p}\leq\\
& &\|d \widetilde{\varphi}(x)-d\widetilde{\varphi}_{n}(x)\|_{p}+
\|d\widetilde{\varphi}_{n}(x)-d\widetilde{\varphi}_{n}(0)\|_{p}+
\|d\widetilde{\varphi}_{n}(0)-d\widetilde{\varphi}(0)\|_{p}\leq\\
& &2\|d\varphi-d\varphi_{n}\|_{\infty}+
\|d\widetilde{\varphi}_{n}(x)-d\widetilde{\varphi}_{n}(0)\|_{p}+
2\|d\varphi-d\varphi_{n}\|_{\infty}=\\ &
&4\|d\varphi-d\varphi_{n}\|_{\infty}+
\|d\widetilde{\varphi}_{n}(x)-d\widetilde{\varphi}_{n}(0)\|_{p}
\hspace{5.7cm} (2)
\end{eqnarray*}
for all $n\in\mathbb{N}$, $x\in B(0_{p},r)\subset TM_{p}$. Since
$\|d\varphi-d\varphi_{n}\|_{\infty}\to 0$ as $n\to\infty$ we can
find $n_{0}\in\mathbb{N}$ so that
$$\|d\varphi-d\varphi_{n_{0}}\|_{\infty}\leq\varepsilon/8.
\eqno(3)$$ Finally, since $d\widetilde{\varphi}_{n_{0}}$ is
continuous at $0_{p}$, there exists $\delta\in (0,r)$ such that
    $$
    \|d\widetilde{\varphi}_{n_{0}}(x)-
    d\widetilde{\varphi}_{n_{0}}(0)\|_{p}\
    \leq\frac{\varepsilon}{2} \eqno(4)
    $$
if $\|x\|_{p}\leq\delta$. By combining $(2)$, $(3)$ and $(4)$, we
get that
$\|d\widetilde{\varphi}(x)-d\widetilde{\varphi}(0)\|_{p}\leq\varepsilon$
if $\|x\|_{p}\leq\delta$. This shows that $d\widetilde{\varphi}$
is continuous at $0_{p}$.
\end{proof}

\medskip

Now the proof of Theorem \ref{Smooth Variational Principle} is an
obvious combination of the above Lemmas.

\begin{rem}
{\em It should be noted that Lemma \ref{the superlemma} is quite a
powerful statement from which a lot of other perturbed
minimization principles can be obtained. For instance:
\begin{enumerate}
\item When we take $M=X$, a complete metric space, and $Y$ is the
space of all the Lipschitz and bounded functions
$f:X\to\mathbb{R}$, with the norm
    $$
    \|f\|_{Y}=\|f\|_{\infty}+\textrm{Lip}(f)=
    \|f\|_{\infty}+\sup\{\frac{|f(x)-f(y)|}{d(x,y)} : x,y\in X,
    x\neq y\}
    $$
(which satisfies $(1)$ and $(2)$ of Lemma \ref{the superlemma}
with $C=1$ and any $r$), then we obtain a statement that is easily
seen to imply Ekeland's variational principle.
\item  When we consider $M=X$, a Banach space having a $C^1$
smooth Lipschitz bump function, and we define $Y$ as the Banach
space of $C^1$ smooth Lipschitz functions $f:X\to\mathbb{R}$, with
the norm
    $$
    \|f\|_{Y}=\|f\|_{\infty}+\|f'\|_{\infty},
    $$
then we recover the known DGZ smooth variational principle for
Banach spaces.
\item Let $M=X$ be any metric space in which some notion of
{\em differentiability} has been defined, and $Y$ be a Banach
space of {\em differentiable} (whatever this word should mean in
this context) and Lipschitz functions $f:X\to\mathbb{R}$, with the
norm
    $$
    \|f\|_{Y}=\|f\|_{\infty}+\textrm{Lip}(f).
    $$
Suppose that $X$ is uniformly bumpable in the sense that $Y$
satisfies $(2)$ of Lemma \ref{the superlemma}. Then we get a
perturbed minimization principle with functions which are {\em
differentiable} and Lipschitz.
\end{enumerate}
}
\end{rem}
\begin{op}
{\em Is Theorem \ref{Smooth Variational Principle} true if one
drops the assumption that $M$ is uniformly bumpable?}
\end{op}

\medskip


\section{A notion of viscosity subdifferential for functions defined on
Riemannian manifolds}

\medskip

\noindent {\bf Definitions and basic properties}

\begin{defn}
{\em Let $(M,g)$ be a Riemannian manifold, and $f:M\longrightarrow
(-\infty,\infty]$ be a proper function. We will say that $f$ is
subdifferentiable at a point $p\in\text{dom}(f)=\{x\in M :
f(x)<\infty\}$ provided there exists a $C^1$ function
$\varphi:M\longrightarrow\mathbb{R}$ such that $f-\varphi$
attains a local minimum at the point $p$. In this case we will
say that $\zeta=d\varphi(p)\in(TM_{p})^{*}\simeq H^{*}=H$ is a
subdifferential of $f$ at $p$. We define the subdifferential set
of $f$ at $p$ by
    $$
    D^{-}f(p)=\{d\varphi(p): \varphi\in C^{1}(M,\mathbb{R}), f-\varphi
    \text{ attains a local minimum at } p\},
    $$
a subset of $T^{*}M_{p}$. Similarly, we define
    $$
    D^{+}f(p)=\{d\varphi(p): \varphi\in C^{1}(M,\mathbb{R}), f-\varphi
    \text{ attains a local maximum at } p\},
    $$
and we say that $f$ is superdifferentiable at $p$ provided
$D^{+}f(p)\neq\emptyset$.

For every $\zeta\in D^{-}f(p)\cup D^{+}f(p)$, the norm of $\zeta$
is defined as
    $$
    \|\zeta\|_{p}=\sup\{|\zeta(h)| : h\in TM_{p}, \|h\|_{p}=1\}.
    $$}
\end{defn}

\begin{rem}\label{trivial remarks on definition of subdifferential}
{\em The following properties are obvious from the definition:
\begin{enumerate}
\item $f$ is subdifferentiable at $p$ if and only if
$-f$ is superdifferentiable at $p$, and
    $$
    D^{+}(-f)(p) = -D^{-}f(p).
    $$
\item If $f$ has a local minimum at $p$ then $0\in D^{-}f(p)$.
\item If $h$ has a local maximum at $p$ then $0\in D^{+}f(p)$.
\end{enumerate}
}
\end{rem}

Next we give other useful equivalent definitions of
subdifferentiability.

\begin{thm}[Characterizations of subdifferentiability]
\label{equivalent definitions of subdifferential} Let $f:M\to
(-\infty, \infty]$ be a function defined on a Riemannian
manifold, $p\in M$, and $\eta\in T^{*}M_{p}$. The following
statements are equivalent:
\begin{enumerate}
\item $\eta\in D^{-}f(p)$, that is, there exists a $C^1$ smooth function
$\varphi:M\to\mathbb{R}$ so that $f-\varphi$ attains a local
minimum at $p$, and $\eta=d\varphi(p)$.
\item There exists a function
$\varphi:M\to\mathbb{R}$ so that $f-\varphi$ attains a local
minimum at $p$, $\varphi$ is Fr{\'e}chet differentiable at $p$,
and $\eta=d\varphi(p)$.
\item For every chart $h:U\subset M\to H$ with $p\in U$, if we take
$\zeta=\eta\circ dh^{-1}(h(p))$ then we have that
    $$
    \liminf_{v\to 0}\frac{(f\circ h^{-1})(h(p)+v)-f(p)-\langle\zeta,v\rangle}{\|v\|}
    \geq 0.
    $$
\item There exists a chart $h:U\subset M\to H$ with $p\in U$ and
such that, for $\zeta=\eta\circ dh^{-1}(h(p))$, we have
    $$
    \liminf_{v\to 0}\frac{(f\circ h^{-1})(h(p)+v)-f(p)-\langle\zeta,v\rangle}{\|v\|}
    \geq 0.
    $$
\end{enumerate}
Moreover, if the function $f$ is locally bounded below (that is,
for every $x\in M$ there is a neighborhood $U$ of $x$ such that
$f$ is bounded below on $U$), then the above conditions are also
equivalent to the following one:
\begin{enumerate}
\item[{(5)}] There exists a $C^1$ smooth function
$\varphi:M\to\mathbb{R}$ so that $f-\varphi$ attains a {\em
global} minimum at $p$, and $\eta=d\varphi(p)$.
\end{enumerate}

\noindent {\em Consequently, any of these statements can be taken
as a definition of $\eta\in D^{-}f(p)$.} Analogous statements are
equivalent in the case of a superdifferentiable function; in
particular $\zeta\in D^{+}f(p)$ if and only if there exists a
chart $h:U\subset M\to H$ with $p\in U$ and such that, for
$\zeta=\eta\circ dh^{-1}(h(p))$,
    $$
    \limsup_{v\to 0}\frac{(f\circ h^{-1})(h(p)+v)-f(p)-\langle\zeta,v\rangle}{\|v\|}
    \leq 0.
    $$
\end{thm}
\begin{proof}

\noindent $(1)\implies (2)$ and $(3)\implies (4)$ are obvious.

\noindent $(2)\implies (3)$. If $f-\varphi$ has a local minimum
at $p$ then $g:=f\circ h^{-1}-\varphi\circ h^{-1}$ has also a
local minimum at $h(p)$, which implies
    $$
    \liminf_{v\to 0}\frac{g(h(p)+v)-g(h(p))}{\|v\|}
    \geq 0
    $$
and, by combining this inequality with the fact that
    $$
    \lim_{v\to 0}\frac{(\varphi\circ h^{-1})(h(p)+v)-
    (\varphi\circ h^{-1})(h(p))-\langle\zeta, v\rangle}{\|v\|}=0
    $$
(because $\zeta=d(\varphi\circ h^{-1})(h(p))$), it is easily
deduced that
    $$
    \liminf_{v\to 0}\frac{(f\circ h^{-1})(h(p)+v)-(f\circ h^{-1})(h(p))-
    \langle\zeta,v\rangle}{\|v\|}
    \geq 0.
    $$

\noindent $(4)\implies (1)$. In order to prove this we will use
the following lemma, which is shown in \cite{DGZ} in a more
general situation.
\begin{lem}
If $V$ is an open set of a Hilbert space $H$, $x\in V$, and
$F:V\to (-\infty, \infty]$ is a function satisfying
    $$
    \liminf_{v\to 0}\frac{F(x+v)-F(x)-\langle \tau, v\rangle}{\|v\|}\geq 0
    $$
for some $\tau\in H^{*}$, then there exists a $C^1$ smooth
function $\psi:H\to\mathbb{R}$ such that $F-\psi$ has a local
minimum at $x$, and $d\psi(x)=\tau$.
\end{lem}
Take an open neighborhood $V$ of $p$ so that $\overline{V}\subset
U$. Note that $F:=f\circ h^{-1}$ is a function from the open
subset $h(U)$ of the Hilbert space $H$ into $(-\infty, \infty]$,
and by the hypothesis we have that
    $$
    \liminf_{v\to
    0}\frac{F(h(p)+v)-F(h(p))-\langle\zeta,v\rangle}{\|v\|}
    \geq 0.
    $$
By the preceding lemma, there exists a $C^{1}$ smooth function
$\psi:h(U)\to\mathbb{R}$ such that $F-\psi$ has a local minimum
at $h(p)$ and $\zeta=d\psi(h(p))$. Let us define $\phi:=\psi\circ
h:U\to\mathbb{R}$, which is a $C^1$ smooth function. It is clear
that $F\circ h-\psi\circ h=f-\phi$ has a local minimum at $p$, and
$d\phi(p)=d\psi(h(p))\circ dh(p)=\zeta\circ dh(p)=\eta$. In order
to finish the proof it is enough to extend $\phi$ to the
complement of $V$ by defining $\varphi=\theta\phi$, where $\theta$
is a $C^1$ smooth Uryshon-type function which is valued $1$ on the
set $V$ and $0$ outside $U$ (such a function certainly exists
because $M$ has $C^\infty$ smooth partitions of unity and
$\overline{V}\subset U)$. It is obvious that $\varphi$ keeps the
relevant properties of $\phi$.

Finally, let us see that, when $f$ is locally bounded below,
$(1)\iff (5)$. Obviously, $(5)\implies (1)$. To see that
$(1)\implies (5)$, let us assume that there exists a $C^1$ smooth
function $\psi:M\to\mathbb{R}$ and some $r>0$ such that
$0=f(p)-\psi(p)\leq f(x)-\psi(x)$ if $x\in B(p,r)$, and denote
$\eta=d\psi(p)$. We have to see that there exists a $C^1$ smooth
function $\varphi:M\to\mathbb{R}$ such that $f-\varphi$ attains a
{\em global} minimum at $p$ and $d\varphi(p)=\eta$.

Consider the open set $U=M\setminus\overline{B}(p,r/2)$.  Since
$f-\psi$ is locally bounded below, for each $x\in U$ there exist
$\delta_{x}>0$ and $m_{x}\in\mathbb{R}$ such that
$B(x,\delta_{x})\subset U$ and $m_{x}\leq f(y)-\psi(y)$ for all
$y\in B(x,\delta_{x})$. Consider the open covering
$$G:=\{B(x,\delta_{x}) : x\in U\}\cup \{B(p,r)\}$$ of $M$. Since
$M$ has $C^\infty$ smooth partitions of unity there exists a
locally finite refinement $\{U_{i}\}_{i\in I}$ of the covering
$G$ and a family of functions $\{\psi_{i}\}_{i\in I}\subset
C^{\infty}(M, [0,1])$ so that $\textrm{supp}(\psi_{i})\subset
U_{i}$ for each $i$ and $\sum_{i\in I}\psi_{i}=1$.

For each $i\in I$, if $U_{i}\subset B(p,r)$ then we define
$\alpha_{i}=0$. Otherwise we can pick some $x_{i}\in U=
M\setminus\overline{B}(p,r/2)$ such that $U_{i}\subset
B(x_{i},\delta_{x_{i}})$, and in this case we define
$\alpha_{i}=m_{x_{i}}$. Now we can define our function
$\varphi:M\to\mathbb{R}$ by
    $$
    \varphi(x)=\psi(x)+\sum_{i\in I}\alpha_{i}\psi_{i}(x).
    $$
It is clear that $\varphi$ is a $C^1$ smooth function such that
$\varphi=\psi$ on $\overline{B}(p,r/2)$ (indeed, take $x\in
\overline{B}(p,r/2)$; if $x\in U_{i}$ then $U_{i}\subset B(p,r)$
because of the choice of the covering $G$ and the $\delta_{y}$, so
$\alpha_{i}=0$, while for all the rest of $j\in I$ we have
$\psi_{j}(x)=0$; therefore $\varphi(x)=\psi(x)+0=\psi(x)$). In
particular, it follows that $\eta=d\psi(p)=d\varphi(p)$.

We claim that $f-\varphi$ attains a global minimum at $p$. Indeed,
fix $x\in M$. If $x\in\overline{B}(p,r/2)=M\setminus U$ then, as
we have just seen, $\varphi(x)=\psi(x)$, and $0=(f-\varphi)(p)\leq
(f-\psi)(x)=(f-\varphi)(x)$. If $x\in U$ then, for those $i\in I$
such that $x\in U_{i}$ we have $(f-\psi)(x)\geq
m_{x_{i}}=\alpha_{i}$, while $\psi_{j}(x)=0$ for those $j\in I$
with $x\notin U_{j}$. Therefore,
\begin{eqnarray*}
& &f(x)-\varphi(x)=f(x)-\psi(x)-\sum_{i\in
I}\alpha_{i}\psi_{i}(x)=\\ &
&f(x)-\psi(x)-\sum\{\alpha_{i}\psi_{i}(x) : i\in I, x\in
U_{i}\}\geq\\ & &\sup\{\alpha_{i} : i\in I, x\in U_{i}\}-
\sum\{\alpha_{i}\psi_{i}(x) : i\in I, x\in U_{i}\}\geq
0=f(p)-\varphi(p),
\end{eqnarray*}
and $f-\varphi$ has a global minimum at $p$.
\end{proof}

\begin{cor}\label{localization}
Let $f:M\to (-\infty, \infty]$ be a function defined on a
Riemannian manifold, and let $h:U\subset M\to h(U)\subset H$ be a
chart of $M$. Then,
    \begin{eqnarray*}
    & &D^{-}f(p)=\{\zeta\circ dh(p) : \zeta\in H^{*},
    \liminf_{v\to 0}\frac{(f\circ h^{-1})(h(p)+v)-f(p)-\langle\zeta,v\rangle}{\|v\|}
    \geq 0\}\\
    & &=\{\zeta\circ dh(p) : \zeta\in D^{-}(f\circ h^{-1})(h(p))\}.
    \end{eqnarray*}
\end{cor}

\noindent Now we can show that subdifferentiable plus
superdifferentiable equals differentiable.

\begin{prop}\label{subdifferentiable plus superdifferentiable equals
differentiable} A function $f$ is differentiable at $p$ if and
only if $f$ is both subdifferentiable and superdifferentiable at
$p$. In this case, $\{d f(p)\}=D^{-}f(p)=D^{+}f(p)$.
\end{prop}
\begin{proof}
Assume first that $f$ is both subdifferentiable and
superdifferentiable at $p$. Then there exist $C^1$ functions
$\varphi, \psi:M\to\mathbb{R}$ such that $f-\varphi$ and $f-\psi$
have a local minimum and a local maximum at $p$, respectively. We
can obviously assume $f(p)=\varphi(p)=\psi(p)$. Then these
conditions mean that $f(x)-\varphi(x)\geq 0$ and $f(x)-\psi(x)\leq
0$ for all $x\in U$, where $U$ is an open neighborhood of $p$. On
the other hand, $(f-\varphi)-(f-\psi)=\psi-\varphi$ has a local
minimum at $p$, hence $0=d(\psi-\varphi)(p)=d\psi(p)-d\varphi(p)$.
That is, we have that
    $$
    \varphi(x)\leq f(x)\leq \psi(x) \textrm{ for all } x\in U, \,
    \varphi(p)=\psi(p)=f(p), \, \textrm{ and } \, d\varphi(p)=d\psi(p).
    $$
By using charts, it is an easy exercise to check that these
conditions imply that $f$ is differentiable at $p$, with
$df(p)=d\psi(p)=d\varphi(p)$; in particular this argument shows
that $\{d f(p)\}=D^{-}f(p)=D^{+}f(p)$.

Now, if $f$ is differentiable at $p$ then, by the chain rule, so
is $f\circ h^{-1}$ at $h(p)$ for any chart $h:U\subset M\to H$; in
particular, putting $\zeta=d(f\circ h^{-1})(h(p))$, we have
$$\lim_{v\to 0}\frac{(f\circ
h^{-1})(h(p)+v)-f(p)-\langle\zeta,v\rangle}{\|v\|}=0, $$ which,
thanks to Theorem \ref{equivalent definitions of subdifferential},
yields $df(p)=\zeta\circ dh(p)\in D^{-}f(p)\cap D^{+}f(p)$.
\end{proof}

What the above proof really shows is the (not completely obvious)
following result: a function $f$ is differentiable at a point $p$
if and only if its graph is trapped between the graphs of two
$C^1$ smooth functions which have the same derivative at $p$ and
touch the graph of $f$ at $p$.
\begin{cor}[Criterion for differentiability]
\label{characterization of differentiability} A function
$f:M\to\mathbb{R}$ is Fr{\'e}chet differentiable at a point $p$
if and only if there are $C^{1}$ smooth functions $\varphi,
\psi:M\to\mathbb{R}$ such that
    $$
    \varphi(x)\leq f(x)\leq \psi(x) \textrm{ for all } x\in M, \,
    \varphi(p)=\psi(p)=f(p), \, \textrm{ and } \, d\varphi(p)=d\psi(p).
    $$
\end{cor}

Let us say a few words about the relationship between
subdifferentiability and continuity. In general, a
subdifferentiable function need not be continuous. For instance,
the function $f:\mathbb{R}\longrightarrow\mathbb{R}$ defined by
$f(x)=0$ if $x\in [0,1]$, and $1$ elsewhere, is Fr{\'e}chet
subdifferentiable everywhere in $\mathbb{R}$, and yet $f$ is not
continuous at $0$ and $1$. However, it is easy to see that
subdifferentiability implies lower semicontinuity.
\begin{prop}\label{subdifferentiable implies lower semicontinuous}
If $f$ is subdifferentiable at $p$ then $f$ is lower
semicontinuous at $p$. In the same way, superdifferentiability
implies upper semicontinuity.
\end{prop}
\begin{proof}
The result is immediate in the case of a function $g:V\subset
H\to (-\infty, \infty]$; indeed, if
    $$
    \liminf_{v\to 0}\frac{g(x+v)-g(x)-\langle\tau, v\rangle}{\|v\|}\geq 0
    $$
then $\liminf_{y\to x}g(y)\geq g(x)$. The general case follows by
applying Theorem \ref{equivalent definitions of subdifferential}.
\end{proof}

\medskip

\noindent {\bf Some rules and fuzzy rules}

Next we study some properties of the subdifferentials related to
composition, sum and product of subdifferentiable and
differentiable functions. Of course, all the statements hold for
superdifferentials as well, with obvious modifications.

\begin{prop}[Chain rule]\label{chain rule}
Let $M$, $N$ be Riemannian manifolds, $g:M\to N$, and $f:N\to
(-\infty,\infty]$. Assume that the function $f$ is
subdifferentiable at $g(p)$, and that $g$ is Fr{\'e}chet
differentiable at $p$. Then the composition $f\circ g:M\to
(-\infty,\infty]$ is subdifferentiable at $p$, and
    $$
    \{\zeta\circ dg(p) : \zeta\in D^{-}f(g(p))\}\subseteq D^{-}(f\circ g)(p).
    $$
\end{prop}
\begin{proof}
Take $\zeta\in D^{-}f(g(p))$, then there exists a function
$\varphi:N\to\mathbb{R}$ so that $f-\varphi$ has a local minimum
at $g(p)$, $\varphi$ is Fr{\'e}chet differentiable at $g(p)$, and
$\zeta=d\varphi(g(p))$. In particular there exists
$\varepsilon>0$ such that $f(y)-\varphi(y)\geq
f(g(p))-\varphi(g(p))$ whenever $d(y, g(p))<\varepsilon$. Define
$\psi=\varphi\circ g$. Since $g$ is differentiable at $p$ and
$\varphi$ is differentiable at $g(p)$, by the chain rule it
follows that $\psi$ is a function from $M$ into $\mathbb{R}$ which
is Fr{\'e}chet differentiable at $p$, with
$d\psi(p)=d\varphi(g(p))\circ dg(p)$. Since $g$ is continuous at
$p$, there exists $\delta>0$ such that $d(g(x),g(p))<\varepsilon$
for all $x$ with $d(x,p)<\delta$. Then we get
$f(g(x))-\varphi(g(x))\geq f(g(p))-\varphi(g(p))$ if
$d(x,p)<\delta$, that is, $f\circ g-\psi$ has a local minimum at
$p$. By Theorem \ref{equivalent definitions of
subdifferential}$[(1)\iff (2)]$, this ensures that $f\circ g$ is
subdifferentiable at $p$, with $\zeta\circ
dg(p)=d\varphi(g(p))\circ dg(p)=d\psi(p)\in D^{-}(f\circ g)(p).$
\end{proof}
The following example shows that the inclusion provided by
Proposition \ref{chain rule} is strict, in general.
\begin{ex}
{\em Let $M=N=\mathbb{R}$, $g(x)=|x|^{3/2}$, $f(y)=|y|^{1/2}$;
$f\circ g(x)=|x|^{3/4}$. Then $g$ is $C^1$ smooth on $\mathbb{R}$,
and we have $dg(0)=0$, $D^{-}f(g(0))=D^{-}f(0)=(-\infty,\infty)$,
$D^{-}(f\circ g)(0)=(-\infty,\infty)$. Therefore $\zeta\circ
dg(0)=0$ for every $\zeta\in D^{-}f(g(0))$.}
\end{ex}

\begin{cor}\label{preservation of subdifferentials by diffeomorphisms}
Let $M$, $N$ be Riemannian manifolds, $h:M\to N$ a $C^1$
diffeomorphism. Then, $f:M\to (-\infty,\infty]$ is
subdifferentiable at $p$ if and only if $f\circ h^{-1}$ is
subdifferentiable at $h(p)$, and
    $$
    D^{-}f(p)=\{\zeta\circ dh(p) : \zeta\in D^{-}(f\circ h^{-1})(h(p))\}.
    $$
\end{cor}
\begin{proof}
If $f:M\to (-\infty,\infty]$ is subdifferentiable at $p$ then, by
the preceding Proposition, $f\circ h^{-1}:N\to (-\infty,\infty]$
is subdifferentiable at $h(p)\in N$ and, moreover, we know that if
$T\in D^{-}f(p)$ then $\zeta:=T\circ dh^{-1}(h(p))\in
D^{-}(f\circ h^{-1})(h(p))$. Then $T=\zeta\circ dh(p)$, with
$\zeta\in D^{-}(f\circ h^{-1})(h(p))$.

Conversely, if $f\circ h^{-1}$ is subdifferentiable at $h(p)$
then, again by the preceding result, $f=(f\circ h^{-1})\circ h$
is subdifferentiable at $p$ and, for any $\zeta\in D^{-}(f\circ
h^{-1})(h(p))$, we have $\zeta\circ dh(p)\in D^{-1}\big((f\circ
h^{-1})\circ h\big)(p)=D^{-}f(p)$.
\end{proof}

\begin{prop}[Sum rule]\label{sum rule}
For all functions $f_{1}, f_{2}:M\longrightarrow
(-\infty,\infty]$, $p\in M$, we have
    $$
    D^{-}f_{1}(p)+D^{-}f_{2}(p)\subseteq D^{-}(f_{1}+f_{2})(p),
    $$
and analogous inclusions hold for superdifferentials.
\end{prop}
\begin{proof}
Take $\zeta_{i}\in D^{-}f_{i}(p)$, $i=1,2$. There are $C^1$
smooth functions $\varphi_{i}:M\to\mathbb{R}$ such that
$f_{i}-\varphi_{i}$ have a minimum at $p$ and
$\zeta_{i}=d\varphi_{i}(p)$ for $i=1,2$. Then
$(f_{1}+f_{2})-(\varphi_{1}+\varphi_{2})=(f_{1}-\varphi_{1})+
(f_{2}-\varphi_{2})$ clearly has a minimum at $p$, hence
$\zeta_{1}+\zeta_{2}=d(\varphi_{1}+\varphi_{2})(p)$ belongs to
$D^{-}(f_{1}+f_{2})(p)$.
\end{proof}
When one of the functions involved in the sum is uniformly
continuous the inclusion provided by this statement can be
reversed in a fuzzy way. This assumption is necessary in general,
as a counterexample (in the Hilbert space) of Deville and Ivanov
shows; see \cite{DeIv}.
\begin{thm}[Fuzzy rule for the subdifferential of the sum]
\label{Fuzzy rule for the subdifferential of the sum}
 Let $(M,g)$ be a Riemannian manifold. Let $f_{1},
f_{2}:M\longrightarrow\mathbb{R}$ be such that $f_{1}$ is lower
semicontinuous and $f_{2}$ is uniformly continuous. Take $p\in M$,
a chart $(U,\varphi )$ with $p\in U$, $\zeta\in
D^{-}(f_{1}+f_{2})(p)$, $\varepsilon>0$, and a neighborhood $V$ of
$(p,\zeta)$ in the cotangent bundle $T^{*}M$. Then there exist
$p_{1}, p_{2}\in U$, $\zeta_{1}\in D^{-}f_{1}(p_{1})$,
$\zeta_{2}\in D^{-}f_{2}(p_{2})$ such that:
$(p_{i},\zeta_{1}\circ d\varphi (p_1)^{-1}\circ d\varphi (p_i)
+\zeta_{2}\circ d\varphi (p_2)^{-1}\circ d\varphi (p_i))\in V$
for $i=1,2$; and
 $|f_{i}(p_{i})-f_{i}(p)|<\varepsilon$ for $i=1,2$.
\end{thm}
\begin{proof}
Fix a chart $(U,\varphi )$ such that $p\in U$ and $T^*U$\ is
diffeomorphic to $U\times H^*$ through the canonical
diffeomorphism $L:T^*U \to U\times H^*$\ defined by $L(q,\xi
)=(q, \xi \circ d\varphi (q)^{-1})$. The theorem can be
reformulated as follows: {\it for every $p\in U$, $\zeta\in
D^{-}(f_{1}+f_{2})(p)$, and $\varepsilon>0$, there exist $p_{1},
p_{2}\in U$, $\zeta_{1}\in D^{-}f_{1}(p_{1})$, $\zeta_{2}\in
D^{-}f_{2}(p_{2})$ such that: $d(p_1,p_2)<\varepsilon$, $||
\zeta_{1}\circ d\varphi (p_1)^{-1}+\zeta_{2}\circ d\varphi
(p_2)^{-1}-\zeta \circ d\varphi (p)^{-1}||<\varepsilon$, and
$|f_{i}(p_{i})-f_{i}(p)|<\varepsilon$ for $i=1,2$.} But this
statement follows immediately from Deville and El Haddad's fuzzy
rule for Banach spaces \cite{DH1} applied to the functions
$f_1\circ \varphi ^{-1}$\ and $f_2\circ \varphi ^{-1}$.
\end{proof}

\begin{prop}[Product rule]\label{product rule}
Suppose $f_{1}, f_{2}:M\longrightarrow [0,\infty)$ are functions
subdifferentiable at $p\in M$. Then $f_{1} f_{2}$ is
subdifferentiable at $p$, and
    $$
    f_{1}(p)D^{-}f_{2}(p)+f_{2}(p)D^{-}f_{1}(p)\subseteq D^{-}(f_{1}f_{2})(p).
    $$
\end{prop}
\begin{proof}
If one of the functions vanishes at $p$ then it has a minimum at
$p$, and so does its product with the other function, hence the
result obviously holds true (recall that if $g$ has a minimum at
$p$ then $0\in D^{-}g(p)$). Therefore we may assume that
$f_{i}(p)>0$ for $i=1,2$.

Pick $\zeta_{i}\in D^{-}f_{i}(p)$, and find $C^1$ smooth functions
$\varphi_{i}:M\to\mathbb{R}$ such that $f_{i}-\varphi_{i}$ have a
local minimum at $p$ and $\zeta_{i}=d\varphi_{i}(p)$ for $i=1,2$.
As usual we may assume that $\varphi_{i}(p)=f_{i}(p)$, so that
$f_{i}-\varphi_{i}\geq 0$. Since $\varphi_{i}(p)=f_{i}(p)>0$ and
$\varphi_{i}$ is continuous, there exists a neighborhood $V$ of
$p$ such that $\varphi_{i}\geq 0$ on $V$. We may assume that $V$
is small enough so that the restrictions of $\varphi_{i}-f_{i}$ to
$V$ have a global minimum at $p$. Then we deduce that
$f_{1}f_{2}\geq \varphi_{1}f_{2}\geq \varphi_{1}\varphi_{2}\,
    \textrm{ on } \, V,$
that is, $$(f_{1}f_{2}-\varphi_{1}\varphi_{2})(x)\geq 0=
(f_{1}f_{2}-\varphi_{1}\varphi_{2})(p) \,
    \textrm{ for all } \, x\in V,$$ which means that
$f_{1}f_{2}-\varphi_{1}\varphi_{2}$ has a local minimum at $p$,
and therefore
    $$
    f_{1}(p)\zeta_{2}+f_{2}(p)\zeta_{1}=
    \varphi_{1}(p)d\varphi_{2}(p)+\varphi_{2}(p)d\varphi_{1}(p)=
    d(\varphi_{1}\varphi_{2})(p)\in D^{-}(f_{1}f_{2})(p).
    $$
\end{proof}
\begin{rem}
{\em If the functions are not positive, the result is not
necessarily true, as the following example shows: $M=\mathbb{R}$,
$f_{1}(x)=|x|$, $f_{2}(x)=-1$, $p=0$ (note that the function
$(f_{1}f_{2})(x)=-|x|$ is not subdifferentiable at $0$).}
\end{rem}

\medskip

\noindent {\bf Topological and geometrical properties of the
subdifferential sets}.

\begin{prop}\label{the subdifferential set is closed and convex}
$D^{-}f(p)$ and $D^{+}f(p)$ are closed and convex subsets of
$T^{*}M_{p}$. In particular, if $f$ is locally Lipschitz then
these sets are $w^{*}$-compact as well.
\end{prop}
\begin{proof}
Let us first check that $D^{-}f(p)$ is convex. Pick $\zeta_{1},
\zeta_{2}\in D^{-}f(p)$, and find $C^1$ smooth functions
$\varphi_{1}, \varphi_{2}:M\to\mathbb{R}$ such that
$d\varphi_{i}(p)=\zeta_{i}$, and $(f-\varphi_{i})(x)\geq
0=(f-\varphi_{i})(p)$ for all $x$ in a neighborhood of $p$. Take
$t\in [0,1]$, and define the function
$\varphi_{t}:M\to\mathbb{R}$ by
$\varphi_{t}(x)=(1-t)\varphi_{1}(x)+t\varphi_{2}(x)$. It is
immediately seen that $\varphi_{t}$ is a $C^1$ smooth function
such that $f-\varphi_{t}$ attains a local minimum at $p$, and
therefore
    $$
    (1-t)\zeta_{1}+t\zeta_{2}=d\varphi_{t}(p)\in D^{-}f(p).
    $$

Now let us see that $D^{-}f(p)$ is closed. Take a chart
$h:U\subset M\to H$ with $p\in U$. Since $dh(p):TM_{p}\to H$ is a
linear isomorphism and $(dh(p))^{*}:H^{*}\to (TM_{p})^{*}$
(defined by $(dh(p))^{*}(\zeta)=\zeta\circ dh(p)$) is a linear
isomorphism as well, and, by Corollary \ref{localization}, we
know that $$D^{-}f(p)=\{\zeta\circ dh(p) : \zeta\in D^{-}(f\circ
h^{-1})(h(p))\}=(dh(p))^{*}\big(D^{-}(f\circ h^{-1})(h(p))\big),$$
it is enough to show that $D^{-}(f\circ h^{-1})(h(p))$ is closed
in $H^{*}$. That is, we have to show the result in the case of a
function $g:V\subset H\to (-\infty,\infty]$ which is
subdifferentiable at a point $x$. So let us prove that
$D^{-}g(x)$ is closed in $(H^{\ast},\|\cdot\|)$. Let
$(p_{n})\subset D^{-}g(x)$ be such that $\|p_{n}-p\|\to 0$, and
let us check that $p\in D^{-}g(x)$. We have $$ \liminf_{v\to
0}\frac{g(x+v)-g(x)-\langle p_{n}, v\rangle}{\| v\|}\geq 0 $$ for
all $n$, and therefore
\begin{eqnarray*}
& &\liminf_{v\to 0}\frac{g(x+v)-g(x)-\langle p,v\rangle}{\|
v\|}=\\ & &=\liminf_{v\to 0}\left[\frac{1}{\|
v\|}\left(g(x+v)-g(x)-\langle p_{n}, v\rangle\right)\right.
+\left.\frac{1}{\| v\|}\langle p_{n}-p, v\rangle\right]\\ &
&\geq\liminf_{v\to 0}\frac{1}{\| v\|}\left(g(x+v)-g(x)-\langle
p_{n}, v\rangle\right)+ \liminf_{v\to 0}\frac{1}{\| v\|}\langle
p_{n}-p, v\rangle\\ & &\geq 0+\liminf_{v\to 0}\frac{1}{\|
v\|}\langle p_{n}-p, v\rangle=-\| p_{n}-p\|\,
\end{eqnarray*}
for all $n$, that is, $$ \liminf_{v\to 0}\frac{g(x+v)-g(x)-\langle
p,v\rangle}{\| v\|}\geq-\| p_{n}-p\| $$ for all $n\in\mathbb{N}$,
and since $\| p_{n}-p\|\to 0$ we deduce that $$ \liminf_{v\to
0}\frac{g(x+v)-g(x)-\langle p,v\rangle}{\| v\|}\geq 0, $$ which
means $p\in D^{-}g(x)$.

Finally, when $f$ is locally Lipschitz, by composing with the
inverse of the exponential map (which provides a Lipschitz chart
on a neighborhood of each point) and using Corollary
\ref{localization}, it is easily seen that $D^{-}f(p)$ and
$D^{+}f(p)$ are bounded. Then, by the Alaoglu-Bourbaki Theorem it
follows that these sets are $w^{*}$-compact.
\end{proof}

\medskip

\noindent {\bf Density of the points of subdifferentiability}

As a consequence of the Smooth Variational Principle, every lower
semicontinuous function is subdifferentiable on a dense subset of
its domain.
\begin{prop}\label{all functions are subdifferentiable in a dense set}
Let $M$ be a Riemannian manifold. If $f:M\longrightarrow (-\infty,
\infty]$ is lower semicontinuous and proper then
$\{p\in\text{dom}(f) : D^{-}f(p)\neq\emptyset\}$ is dense in
$\text{dom(f)}:=\{x\in M : f(x)<\infty\}$.
\end{prop}
\begin{proof}
Assume first that $M$ is complete and uniformly bumpable (such
is, for instance, the case when $M$ is a Hilbert space $H$). In
this case we can deduce the result directly by applying the Smooth
Variational Principle \ref{Smooth Variational Principle} as
follows. Pick any point $p_{0}$ with $f(p_{0})<\infty$, and any
open neighborhood $U$ of $p_{0}$. We must show that there is some
$p\in U$ such that $D^{-}f(p)\neq\emptyset$. Since $M$ has smooth
partitions of unity, there is a $C^\infty$ smooth function
$b:M\to [0,\infty)$ such that $b(y)>0$ if and only if $y\in U$.
Consider the function $g:M\to (-\infty, \infty]$ defined by
    $$
    g(y)=\frac{1}{b(y)} \, \textrm{ if } y\in U, \, \textrm{ and }
    \, g(y)=\infty \, \textrm{ if } \, y\notin U.
    $$
The function $g$ is lower semicontinuous on $M$, and $C^\infty$
smooth on $U$. Then the sum $f+g$ is lower semicontinuous, and
$(f+g)(p_{0})<+\infty$. According to the Smooth Variational
Principle, there exists a $C^1$ smooth function
$\varphi:M\to\mathbb{R}$ such that $(f+g)-\varphi$ attains a
strong minimum at some point $p\in M$. In fact we have $p\in U$,
because this function is valued $+\infty$ outside $U$. But, since
the function $\varphi-g$ is $C^1$ smooth on $U$, and
$f-(\varphi-g)$ attains its minimum at $p$, we conclude that
    $$
    d(\varphi-g)(p)\in D^{-}f(p)\neq\emptyset.
    $$

Now let us consider the case when $M$ is not necessarily complete
or uniformly bumpable. Pick a point $p_{0}\in\textrm{dom}(f)$ and
an open set $U$ containing $p_{0}$. We may assume that $U$ is
small enough so that there is a chart $h:\overline{U}\subset V\to
H$. By Corollary \ref{localization} we know that, for any $p\in
M$, we have $D^{-}f(p)\neq\emptyset$ if and only if $D^{-}(f\circ
h^{-1})(p)\neq\emptyset$, so it is enough to see that there is
some $x\in h(U)$ with $D^{-}(f\circ h^{-1})(x)\neq\emptyset$.
Define $F(x)=f\circ h^{-1}(x)$ if $x\in h(\overline{U})$, and
$F(x)=+\infty$ otherwise. The function $F$ is lower semicontinuous
on $H$, and $F=f\circ h^{-1}$ on $h(U)$. Since the Hilbert space
$H$ is certainly complete and uniformly bumpable, we can apply the
first part of the argument to find some $x\in h(U)$ so that
$\emptyset\neq D^{-}F(x)=D^{-}(f\circ h^{-1})(x)$.
\end{proof}

\medskip

\noindent {\bf Mean value inequalities}

There are many subdifferential mean value inequality theorems for
functions defined on Banach spaces. Here we will only consider two
of them, which complement each other. The first one is due to
Deville \cite{DeTVM} and holds for all lower semicontinuous
functions $f$ defined on an open convex set of a Banach space,
even if they are not required to be everywhere subdifferentiable,
but it demands a bound for {\em all} of the subgradients of the
function at all the points where it is subdifferentiable. The
second one is due to Godefroy (who improved a similar previous
result of Azagra and Deville), see \cite{GoTVM, AD1}, and only
demands the existence of {\em one} subdifferential or
superdifferential which is bounded (by the same constant) at each
point, but it requires the function to satisfy $D^{-}f(x)\cup
D^{+}f(x)\neq\emptyset$ for {\em all} the points $x$ in the
domain of $f$ (an open convex subset of a Banach space).

Next we extend these mean value inequality theorems to the setting
of Riemannian manifolds. The main ideas of the proofs of these
results could be adapted to obtain direct proofs which would be
valid for the case of manifolds, but for shortness we choose here
to deduce them from the Hilbert space case.

\begin{thm}[Deville's mean value inequality]\label{TVM Deville}
Let $(M,g)$ be a Riemannian manifold, and let
$f:M\longrightarrow\mathbb{R}$ be a lower semicontinuous function.
Assume that there exists a constant $K>0$ such that
$\|\zeta\|_{p}\leq K$ for all $\zeta\in D^{-}f(p)$ and $p\in M$.
Then, $$|f(p)-f(q)|\leq K d_{M}(p,q)\, \hbox { for all } p,q\in
M.$$
\end{thm}
\begin{proof}
The result is true in the case when $M=H$ is a Hilbert space
\cite{DeTVM}. For completeness we give a hint of Deville's
argument, which is an instructive application of the Smooth
Variational Principle. By standard arguments it suffices to show
the result locally (see the proof of the general case below). Fix
$x_{0}\in H$. Since $f$ is locally bounded below there are $N,
\delta>0$ so that $f(x)-f(x_{0})\geq -N$ whenever $x\in
B(x_{0},2\delta)$. For fixed $y\in B(x_{0},\delta/4)$,
$\varepsilon>0$, consider the function defined by $F(x)=
f(x)-f(y)-\alpha(\|x-y\|)$ for $\|x-y\|\leq\delta$, and
$F(x)=+\infty$ elsewhere, where $\alpha:[0,\infty)\to [0,\infty)$
is $C^1$ smooth and satisfies $\alpha(t)=(K+\varepsilon)t$ if
$t\leq\delta/2$, $\alpha(\delta)\geq N$, and $\alpha'(t)\geq
K+\varepsilon$ for all $t>0$. If $\inf F<0$, by applying the
Smooth Variational Principle one can get a point $x_{1}\in
B(y,\delta)\setminus\{y\}$ and a subgradient $\zeta\in
D^{-}f(x_{1})$ so that $\|\zeta\|>K$, a contradiction. Hence
$F\geq 0$, and by letting $\varepsilon\to 0$ the local result
follows. See \cite{DeTVM} for the details.

Now consider the general case of a Riemannian  manifold. Fix any
two points $p,q\in M$, and consider a continuous and piecewise
$C^1$ smooth path $\gamma:[0,T]\to M$, parameterized by arc
length, with $\gamma(0)=p, \gamma(T)=q$. Take $\varepsilon>0$.
According to Theorem \ref{properties of exp}, for each $x\in
\gamma([0,T])$ there exists $r_{x}>0$ so that $\exp_{x}:B(0_{x},
2r_{x})\subset TM_{x}\to B(x, 2r_{x})\subset M$ is a $C^\infty$
diffeomorphism so that the derivatives of $\exp_{x}$ and
$\exp_{x}^{-1}$ are bounded by $1+\varepsilon$ on these balls.
Since $\gamma([0,T])$ is compact, there are a finite collection
of points $x_{1}=p, x_{2}, ..., x_{n}=q\in \gamma[0,T]$ so that
    $$
    \gamma\big([0,T]\big)\subset\bigcup_{j=1}^{n}B(x_{j},r_{j}),
    $$
where we denote $r_{j}=r_{x_{j}}$ for short. Set $r=\min\{r_{1},
..., r_{n}\}$, and pick an $m\in\mathbb{N}$ big enough so that
$T/m<r/2$. Define $t_{0}=0<t_{1}=T/m < ... < t_{j}=jT/m < ... <
T=t_{m}$, and consider the points $a_{j},b_{j}$ with
$a_{j}=b_{j-1}=\gamma(t_{j-1})$, $j=1,...,m$, and
$b_{m}=\gamma(t_{m})$.

For each $j\in\{1,..., m-1\}$ we may choose an
$i_{j}\in\{1,...,n\}$ so that $$\gamma[t_{j-1}, t_{j}]\cap
B(x_{i_{j}}, r_{i_{j}})\neq\emptyset,$$ and we also set
$i_{0}=1$, $i_{m}=n$ (so that $x_{i_{0}}=p$ and $x_{i_{m}}=q$).
Since the length of the restriction of $\gamma$ to
$[t_{j-1},t_{j}]$, which we denote $\gamma_{j}$, is
$t_{j}-t_{j-1}=T/m<r/2\leq r_{i_{j}}/2$, this obviously means that
    $$
    \gamma[t_{j-1}, t_{j}]\subset B(x_{i_{j}}, 2r_{i_{j}})
    $$
for each $j=1,...,m$. In order to avoid an unnecessary burden of
notation, in the sequel we denote $y_{j}=x_{i_{j}}$, and
$s_{j}=r_{i_{j}}$, for $j=0, 1, ..., m$.

Consider the function $f_{j}:B(0_{y_{j}}, 2s_{j})\to\mathbb{R}$
defined by $f_{j}=f\circ\exp_{y_{j}}$. By Corollary
\ref{preservation of subdifferentials by diffeomorphisms} we know
that
    $$
    D^{-}f_{j}(x)=
    \{\zeta\circ d\exp_{y_{j}}(x) : \zeta\in D^{-}(f)(\exp_{y_{j}}(x))\}
    $$
for all $x\in B(0_{y_{j}}, 2s_{j})$. Since $\|\zeta\|_{y}\leq K$
for all $\zeta\in D^{-}f(y)$ with $y\in M$, and
$\|d\exp_{y_{j}}(x)\|\leq(1+\varepsilon)$ for all $x\in
B(0_{y_{j}},2s_{j})$, we deduce that $\|\eta\|_{y_{j}}\leq
(1+\varepsilon)K$ for all $\eta\in D^{-}f_{j}(x)$, $x\in
B(0_{y_{j}},2s_{j})$. Then we can apply the result for the case
$H=TM_{x_{j}}$ and the function $f_{j}$ to see that
\begin{eqnarray*}
& &|f(a_{j})-f(b_{j})|=
    |f_{j}(\exp_{y_{j}}^{-1}(a_{j}))-f(\exp_{y_{j}}^{-1}(b_{j}))|\leq
    \\ & &(1+\varepsilon) K d_{TM_{y_{j}}}\big(\exp_{y_{j}}^{-1}(a_{j})),
    \exp_{y_{j}}^{-1}(b_{j})\big)
\end{eqnarray*}
for all $j=1,2,...,m$. On the other hand, since
$\exp_{x_{j}}^{-1}$ is $(1+\varepsilon)$-Lispchitz we also have
    $$
    d_{TM_{y_{j}}}\big(\exp_{y_{j}}^{-1}(a_{j})), \exp_{y_{j}}^{-1}(b_{j})\big)
    \leq (1+\varepsilon)d_{M}(a_{j},b_{j})
    $$
for all $j=1,2,...,m$. By combining these two last inequalities
we deduce that
    $$
    |f(a_{j})-f(b_{j})|\leq (1+\varepsilon)^{2} K d_{M}(a_{j},b_{j})\leq
    (1+\varepsilon)^{2} K \int_{t_{j-1}}^{t_{j}}\|d\gamma(t)\|dt$$
for all $j=1, ..., m$. Therefore,
\begin{eqnarray*}
& &|f(p)-f(q)|=|\sum_{j=1}^{m}(f(a_{j})-f(b_{j}))|\leq
    \sum_{j=1}^{m}|f(a_{j})-f(b_{j})|\leq\\
    & &(1+\varepsilon)^{2} K \sum_{j=1}^{m}\int_{t_{j-1}}^{t_{j}}\|d\gamma(t)\|dt=
    (1+\varepsilon)^{2} K \int_{0}^{T}\|d\gamma(t)\|dt=
    (1+\varepsilon)^{2} K L(\gamma).
\end{eqnarray*}
By taking the infimum over the set of continuous and piecewise
$C^1$ paths $\gamma$ joining $p$ and $q$ with length $L(\gamma)$,
we get
    $$
    |f(q)-f(p)|\leq (1+\varepsilon)^{2} K d_{M}(q,p).
    $$
Finally, by letting $\varepsilon$ go to $0$ we obtain the desired
inequality: $|f(q)-f(p)|\leq K d_{M}(q,p)$.
\end{proof}
\begin{cor}
Let $(M,g)$ be a Riemannian manifold, and let
$f:M\longrightarrow\mathbb{R}$ be a continuous function. Then $$
\sup\left\{\| \zeta\|_{p}: \zeta\in D^{-}f(p), p\in
M\right\}=\sup\left\{\| \zeta\|_{p} : \zeta\in D^{+}f(p), p\in
M\right\}. $$ These quantities are finite if and only if $f$ is
Lipschitz on $M$, and in this case they are equal to the Lipschitz
constant of $f$.
\end{cor}

\begin{thm}[Godefroy's mean value inequality]\label{TVM Godefroy}
Let $(M,g)$ be a Riemannian manifold, and let
$f:M\longrightarrow\mathbb{R}$ be a Borel function such that
    $$
    D^{-}f(p)\cup D^{+}f(p)\neq\emptyset
    $$
for every $p\in M$. Define $\Phi:M\longrightarrow\mathbb{R}$ by
    $$
    \Phi(p)=\inf\{\|\zeta\|_{p} : \zeta\in D^{-}f(p)\cup D^{+}f(p)\}.
    $$
Then, for every path $\gamma:I\to M$ parameterized by arc length,
one has that
    $$
    \mu\big(f(\gamma(I))\big)\leq
    \int_{I}\Phi(\gamma(t))dt.
    $$
\end{thm}
\noindent Here $\mu$ is the Lebesgue measure in $\mathbb{R}$.
\begin{proof}
The result is already proved in the case when $M=H$ is a Hilbert
space, see \cite{GoTVM}. Let us see how the general case can be
deduced. Let us denote $I=[0,T]$. For a given $\varepsilon>0$,
choose points $y_{j}=x_{i_{j}}, a_{j}, b_{j}$, and numbers
$s_{j}=r_{i_{j}}, t_{j}$, exactly as in the proof of Theorem
\ref{TVM Deville}. Let us denote
$f_{j}=f\circ\exp_{y_{j}}:B(0_{y_{j}},2s_{j})\to\mathbb{R}$,
$\gamma_{j}=\exp_{y_{j}}^{-1}\circ
\gamma:I_{j}:=[t_{j-1},t_{j}]\to B(0_{y_{j}}, 2s_{j})\subset
TM_{y_{j}}$, and
    $$
    \Phi_{j}(y)=\inf\{\|\zeta\|_{y_{j}} : \zeta\in D^{-}f_{j}(x)\cup D^{+}f_{j}(x)\}
    $$
for each $x\in B(0_{y_{j}}, 2s_{j})$. Since
$D^{-}f_{j}(x)=\{\zeta\circ d\exp_{y_{j}}(x) : \zeta\in
D^{-}(f)(\exp_{y_{j}}(x))\}$ for all $x\in B(0_{y_{j}},2s_{j})$,
and $\exp_{y_{j}}$ is $(1+\varepsilon)$-bi-Lipschitz on these
balls, it is easy to see that $\Phi_{j}(x)\leq
(1+\varepsilon)\Phi(\exp_{y_{j}}(x))$ for all $x\in
B(0_{y_{j}},2s_{j})$.

By applying the result for $H=TM_{y_{j}}$, the function $f_{j}$
and the path $\gamma_{j}$, we get that
    $$
    \mu\big(f(\gamma(I_{j}))\big)=
    \mu\big(f_{j}(\gamma_{j}(I_{j}))\big)\leq
    \int_{I_{j}}\Phi_{j}(\gamma_{j}(t))dt,
    \eqno(1)
    $$
for all $j=1,2,...,m$. But we also have that
$$\int_{I_{j}}\Phi_{j}(\gamma_{j}(t))dt\leq
\int_{I_{j}}(1+\varepsilon)\Phi(\exp_{y_{j}}(\gamma_{j}(t)))dt=
(1+\varepsilon)\int_{I_{j}}\Phi(\gamma(t))dt. \eqno(2)$$ By
combining inequalities $(1)$ and $(2)$, and summing over
$j=1,...,m$, we get
\begin{eqnarray*}
& &\mu\big(f(\gamma(I))\big)\leq
    \sum_{j=1}^{m}\mu\big(f(\gamma(I_{j}))\big)\leq
    \sum_{j=1}^{m}(1+\varepsilon)
    \int_{I_{j}}\Phi(\gamma(t))dt\leq\\
    & &(1+\varepsilon)\sum_{j=1}^{m}\int_{I_{j}}\Phi(\gamma(t))dt=
    (1+\varepsilon)\int_{I}\Phi(\gamma(t))dt.
\end{eqnarray*}
Finally, by letting $\varepsilon$ go to $0$ we get
$\mu\big(f(\gamma(I))\big)\leq \int_{I}\Phi(\gamma(t))dt.$
\end{proof}

\begin{cor}
Let $(M,g)$ be a Riemannian manifold, $f:M\to\mathbb{R}$ a Borel
function such that for every $p\in M$ there exists $\zeta\in
D^{-}f(p)\cup D^{+}f(p)$ with $\|\zeta\|_{p}\leq K$. Then,
    $$
    \mu\big(f(\gamma(I))\big)\leq K L(\gamma)
    $$
for every path $\gamma:I\to M$. In particular, when $f$ is
continuous it follows that $|f(p)-f(q)|\leq K d_{M}(p,q)$ for all
$p,q\in M$.
\end{cor}

\medskip


\section{(Sub)differentiability of convex functions on Riemannian manifolds}

he aim of this section is to prove that every (continuous) convex
function defined on a Riemannian manifold is everywhere
subdifferentiable, and differentiable on a dense set.

\begin{defn}\label{def convexity by convexity along geodesics}
{\em Let $M$ be a Riemannian manifold. A function
$f:M\rightarrow\mathbb{R}$ is said to be convex provided that the
function $f\circ \sigma:I\subseteq\mathbb{R}\to\mathbb{R}$ is
convex for every geodesic $\sigma:I\rightarrow M$ (parameterized
by arc length).}
\end{defn}

The following Proposition is probably known, at least in the case
when $M$ is finite-dimensional, but we provide a short proof for
the reader's convenience, as we have not been able to find an
explicit reference.
\begin{prop}\label{every convex function is Lipschitz}
Let $M$ be a Riemannian manifold. If a function
$f:M\to\mathbb{R}$ is convex and locally bounded, then $f$ is
locally Lipschitz. In particular, every continuous convex
function is locally Lipschitz.
\end{prop}
\begin{proof}
Take $p\in M$. Since $f$ is locally bounded there exists $R>0$
such that $f$ is bounded on the ball $B(p,R)$. According to
Theorem \ref{Riemannian manifolds are locally convex}, there
exists $r>0$ with $0<r<R/2$ such that the open balls $B(p, 2r)$
and $B(p,r)$ are convex. Fix $C=\sup\{f(x) : x\in B(p, 2r)\}$, and
$m=\inf\{f(x) : x\in B(p, 2r)\}$. We are going to see that $f$ is
$K$-Lispchitz on the ball $B(p,r)$, where $K=(C-m)/r$. Indeed,
take $x_{1}, x_{2}\in B(p,r)$. Since $B(p,r)$ is convex, there
exists a unique geodesic $\gamma:[t_{1}, t_{2}]\to B(p,r)$, with
length $d(x_{1}, x_{2})=t_{2}-t_{1}$, joining $x_{1}$ to $x_{2}$.
Take $v_{1}\in TM_{x_{1}}$ such that
$\gamma(t)=\exp_{x_{1}}\big((t-t_{1})v_{1}\big)$ for $t\geq t_{1}$
small enough. Since the ball $B(p,2r)$ is still convex and
$x_{1}\in B(p,r)$, we may define a geodesic $\sigma_{1}:[-r,r]\to
B(p,2r)\subset M$ through $x_{1}$ by
    $$
    \sigma_{1}(t)=\exp_{x_{1}}(t v_{1}) \textrm{ for all } t\in [-r,r].
    $$
In the same way we may take $v_{2}\in TM_{x_{2}}$ and define a
geodesic $\sigma_{2}:[-r,r]\to B(p,2r)\subset M$ through $x_{2}$
by
    $$
    \sigma_{2}(t)=\exp_{x_{2}}(t v_{2}) \textrm{ for all } t\in
    [-r,r],
    $$
in such a way that
$\gamma(t)=\exp_{x_{2}}\big((t-t_{2})v_{2}\big)$ for $t\leq
t_{2}$ with $|t|$ small enough. Set $t_{3}=t_{1}-r$,
$t_{4}=t_{2}+r$, $x_{3}=\sigma_{1}(-r)$, $x_{4}=\sigma_{2}(r)$,
and $I=[t_{3}, t_{4}]$. Then, if we define $\sigma:I\to B(p, 2r)$
by
    $$
    \sigma(t)=\left\{
    \begin{array}{lll}
    &\sigma_{1}(t-t_{1})\quad \textrm{if $t\in [t_{3}, t_{1}]$;}\\
    &\gamma(t) \quad \hspace{0.9cm} \textrm{if $t\in [t_{1}, t_{2}]$;}\\
    &\sigma_{2}(t-t_{2}) \quad \textrm{if $t\in [t_{2}, t_{4}]$,}
    \end{array}\right.
    $$
it is clear that $\sigma$ is a geodesic joining $x_{3}$ to
$x_{4}$ in $B(p, 2r)$. Now, since $f$ is convex, the function
$g:[t_{3}, t_{4}]\subset\mathbb{R}\to\mathbb{R}$ defined as
    $$
    g(t)=f(\sigma(t))
    $$
is convex. Therefore we have
    $$
    \frac{g(t_{1})-g(t_{3})}{t_{1}-t_{3}}\leq\frac{g(t_{2})-g(t_{1})}{t_{2}-t_{1}}
    \leq\frac{g(t_{4})-g(t_{2})}{t_{4}-t_{2}},
    $$
where $t_{3}=t_{1}-r<t_{1}<t_{2}<t_{2}+r=t_{4}$. Bearing in mind
that $x_{3}, x_{4}\in B(p, 2r)$, and $t_{2}-t_{1}=d(x_{1},
x_{2})$, it follows that
    $$
    -\frac{C-m}{r}\leq \frac{f(x_{1})-f(x_{3})}{r}\leq
    \frac{f(x_{2})-f(x_{1})}{d(x_{1}, x_{2})}
    \leq\frac{f(x_{4})-f(x_{2})}{r}\leq\frac{C-m}{r}.
    $$
This shows that $|f(x_{1})-f(x_{2})|\leq K d(x_{1}, x_{2})$ for
all $x_{1}, x_{2}\in B(p, r)$, where $K=(C-m)/r$.
\end{proof}

Let us recall that for a locally Lipschitz function
$F:H\to\mathbb{R}$ on a Hilbert space $H$, we may define the
generalized directional derivative $F^0(x,v)$ as the
$$\limsup_{(y,t)\to (x,0^{+})}{F(y+tv)-F(y) \over t}.$$ For every
$x\in H$, $F^0(x,v)$ is a subadditive positively homogeneous
function of $v$, and the set $\{ x^*\in H^*: x^*(v)\leq F^0(x,v)\,
\textrm{ for all } v\}$\ is called the generalized gradient of
$F$\ at $x$, and is denoted by $\partial F(x)$. The generalized
gradient is a nonempty, convex, $w^*$-compact subset of $H^*$; see
\cite{C} for more information.

\begin{thm}\label{convex functions are everywhere subdifferentiable}
Let $g:M\to\mathbb{R}$ be a continuous convex function on a
Riemannian manifold. Then $g$ is subdifferentiable at every point
of $M$.
\end{thm}
\begin{proof}
Let $\phi _p:U_p \to H$ be an exponential chart at $p$. We have
$\phi_{p}(p)=0$. Given another point $q\in U_p$, take a $(\phi
_p,v)\in TM_{q}$, and denote $\sigma_{q,v}(t)=\phi _q^{-1} (tw)$,
where $(\phi _p,v)\sim (\phi _q,w)$, which is a geodesic passing
through $q$ with derivative $(\phi _p,v)$. Here, $(\phi _p,v)\sim
(\phi _q,w)$ means that $w=d(\phi _q\circ \phi _p^{-1})(\phi
_p(q))(v)$, or equivalently $v=d(\phi _p\circ \phi
_q^{-1})(0_{q})(w)$.

Let us define $$f^0(p,v)=\limsup_{q\to p \ t \to 0^+}{f(\sigma
_{q,v}(t))-f(q) \over t}$$
\begin{claim}\label{first claim}
{\em We have that $f^0(p,v)=\inf_{t>0}{f(\sigma _{p,v}(t))-f(p)
\over t}$, and consequently $$f^0(p,v)=\inf_{t>0}{(f \circ \phi
_p^{-1}) (tv)-(f \circ \phi _p^{-1}) (0) \over t}.$$}
\end{claim}
\begin{claim}\label{second claim}
{\em There exists $x^*\in H^*$\ such that $x^*(v)\leq f^0(p,v)$\
for every $v\in H$.}
\end{claim}
\noindent From these facts it follows that $$(f \circ \phi
_p^{-1}) (tv)-(f \circ \phi _p^{-1}) (0)-x^*(tv)\geq 0$$ for
every $v\in S_H$, and every $t\in [0,r)$, provided that
$B(0,r)\subset \phi _p(U_p)$. Hence $(f \circ \phi _p^{-1})
-x^*$\ attains a minimum at $0$ and therefore $x^*\in D^-((f
\circ \phi _p^{-1}) )(0)$. We then conclude that
$D^-f(p)\neq\emptyset$\ by Corollary \ref{localization}. This
shows the theorem.
\end{proof}
\noindent {\bf Proof of Claim \ref{first claim}}.  Fix a $\delta
>0$. Since $f\circ\sigma _{q,v}$\ is convex we have that
$$f^0(p,v)=\lim_{\varepsilon \to
0^+}\sup_{d(p,q)\leq\varepsilon\delta}
\sup_{0<t<\varepsilon}{f(\sigma _{q,v}(t))-f(q) \over t}=$$
$$=\lim_{\varepsilon \to 0^+}\sup_{d(p,q)\leq \varepsilon \delta}
{f(\sigma _{q,v}(\varepsilon ))-f(\sigma _{q,v}(0)) \over
\varepsilon}=(*).$$ Next we estimate $d(\sigma _{p,v}(\varepsilon
),\sigma _{q,v}(\varepsilon ))$. We have $$d(\sigma
_{p,v}(\varepsilon ),\sigma _{q,v}(\varepsilon ))\leq K_p||\phi
_p(\sigma _{p,v}(\varepsilon ))-\phi _p(\sigma _{q,v}(\varepsilon
))|| =K_p||\\ \varepsilon v-\phi _p(\sigma _{q,v}(\varepsilon
))||= $$ $$=K_p||\varepsilon v-(\phi _p\circ \phi
_q^{-1})(\varepsilon w)||=K_p||\varepsilon v-(\phi _p\circ \phi
_q^{-1})(0)- \varepsilon d(\phi _p\circ \phi
_q^{-1})(0)(w)-o(\varepsilon )||= $$ $$=K_p||(\phi _p\circ \phi
_q^{-1})(0)+ o(\varepsilon )|| \leq K_p(||\phi
_p(q)||+||o(\varepsilon)||)\leq$$ $$\leq
K_p(L_pd(p,q)+||o(\varepsilon)||)\leq K_p(L_p\varepsilon \delta
+\varepsilon \delta )\leq C\varepsilon \delta,$$ where $L_p$\ and
$K_p$\ are the Lipschitz constants of $\phi _p$\ and $\phi
_p^{-1}$\ respectively, $C=K_{p}(L_{p}+1)$, and $\varepsilon$\ is
small enough so that $\|o(\varepsilon)\|\leq \varepsilon \delta$\
and $||D(\phi _p\circ \phi _q^{-1})(v)-v||<\delta$.

Since $f$ is locally Lipschitz there exists $K>0$ so that $f$ is
$K$-Lispchitz on a neighborhood of $p$ which may be assumed to be
$U_p$. From the above estimates we get that, for
$d(p,q)\leq\varepsilon\delta$, $$ \Big|{f(\sigma
_{q,v}(\varepsilon ))-f(\sigma _{q,v}(0)) \over \varepsilon}-
{f(\sigma _{p,v}(\varepsilon ))-f(\sigma _{p,v}(0)) \over
\varepsilon}\Big|\leq $$ $$\leq {1\over \varepsilon}(|f(\sigma
_{q,v}(\varepsilon ))-f(\sigma
_{p,v}(\varepsilon))|+|f(p)-f(q)|)\leq$$ $$\leq {K\over
\varepsilon}(d(\sigma _{p,v}(\varepsilon ),\sigma
_{q,v}(\varepsilon ))+d(p,q))\leq K(C+1)\delta.$$ Now we deduce
that $$(*)\leq \lim_{\varepsilon \to 0^+}{f(\sigma
_{p,v}(\varepsilon ))-f(\sigma _{p,v}(0))\over
\varepsilon}+K(C+1)\delta$$ and, by letting $\delta \to 0$, we get
$$f^0(p,v)\leq \lim_{t\to 0^+}{f(\sigma _{p,v}(t))-f(p) \over t}=
\inf_{t>0}{f(\sigma _{p,v}(t))-f(p) \over t}.$$ Since the other
inequality holds trivially the claim is proved. \qed

\noindent {\bf Proof of Claim \ref{second claim}}. We have that
$$\limsup_{q\to p \ t\to 0^+}{f(\sigma _{q,v}(t))-f(\sigma
_{q,v}(0)) \over t}= \lim_{\varepsilon \to
0^+}\sup_{d(p,q)<\varepsilon \ 0<t<\varepsilon}{f(\sigma
_{q,v}(t))-f(\sigma _{q,v}(0)) \over t}=$$ $$=\lim_{\varepsilon
\to 0^+}\sup_{d(p,q)<\varepsilon \ 0<t<\varepsilon}{(f \circ \phi
_p^{-1})(\phi _p(\sigma _{q,v}(t)) -(f \circ \phi _p^{-1})(\phi
_p(q)) \over t}=$$ $$=\lim_{\varepsilon \to
0^+}\sup_{d(p,q)<\varepsilon \ 0<t<\varepsilon}{F(y+\lambda
_y(t))-F(y)\over t},$$ where $(f \circ \phi _p^{-1}) =F$, $y=\phi
_p(q)$, and $\lambda _y(t)=\phi _p(\sigma _{q,v}(t))-\phi _p(q)$.
Next we get $$\lim_{\varepsilon \to 0^+}\sup_{d(p,q)<\varepsilon \
0<t<\varepsilon}{F(y+\lambda _y(t))-F(y)\over t}=
   \lim_{\varepsilon \to 0^+}\sup_{||y||<\varepsilon \
   0<t<\varepsilon}{F(y+\lambda _y(t))-F(y)\over t},$$
because $L_p||y||\leq d(p,q)\leq K_p||y||$ (recall that $\phi _p$\
and $(\phi _p)^{-1}$\ are Lipschitz).

Now, if we take $||y||<\varepsilon$\ and $0<t<\varepsilon$, we
have $$\Big|{F(y+\lambda _y(t))-F(y)\over t}-{F(y+tv)-F(y)\over
t}\Big|=
  \Big|{F(y+\lambda _y(t))-F(y+tv)\over t}\Big|\leq $$
$$\leq  K'{||\lambda _y(t)-tv|| \over t} = K'\varphi (t),$$ where
$K'$\ is the Lipschitz constant of $F$ and $\varphi$ satisfies
$\lim_{t\to 0^+}\varphi (t)=0$, because $$\lambda _y(t)-tv=\phi
_p(\sigma _{q,v}(t))-\phi _p(q)-tv=o(t).$$

Finally, we have $$\Big|\sup_{||y||<\varepsilon \
0<t<\varepsilon}{F(y+\lambda _y(t))-F(y)\over t}-
\sup_{||y||<\varepsilon \ 0<t<\varepsilon}{F(y+tv)-F(y)\over
t}\Big|
  \leq $$
$$\leq \sup_{||y||<\varepsilon \ 0<t<\varepsilon}\Big|{F(y+\lambda
_y(t))-F(y)\over t}-{F(y+tv)-F(y)\over t}\Big|\leq$$ $$\leq
K'\sup_{0<t<\varepsilon}\varphi (t),$$ which goes to $0$ as
$\varepsilon \to 0^+$. Therefore $f^0(p,v)=F^0(0,v)$\ and
$x^*(v)\leq f^0(p,v)$\ for every $v\in H$, provided that $x^*\in
\partial F(0)$, the generalized gradient of $F$\ at $0$.  \qed

\begin{thm}
Let $g:M\to\mathbb{R}$ be a continuous convex function on a
complete Riemannian manifold. Then the set Diff$(g):=\{x\in M : g
\textrm{ is differentiable at } x\}$ is dense in $M$.
\end{thm}
\begin{proof}
According to Proposition \ref{all functions are subdifferentiable
in a dense set}, $\textrm{Diff}^{+}(g):=\{p\in M :
D^{+}g(p)\neq\emptyset\}$ is dense in $M$. On the other hand, by
Theorem \ref{convex functions are everywhere subdifferentiable},
we know that $\textrm{Diff}^{-}(g):=\{p\in M :
D^{-}g(p)\neq\emptyset\}=M$. Then, by Proposition
\ref{subdifferentiable plus superdifferentiable equals
differentiable}, we get that
$$\textrm{Diff}(g)=\textrm{Diff}^{+}(g)\cap
\textrm{Diff}^{-}(g)=\textrm{Diff}^{+}(g) \, \textrm{ is dense in
$M$. }$$
\end{proof}

By using more sophisticated tools, this result can be extended to
the category of locally Lipschitz functions, as we next show.

\begin{thm}\label{Lipschitz functions are differentiable somewhere}
Let $g:M\to\mathbb{R}$ be a locally Lipschitz function. If $M$ is
finite-dimensional, then $g$ is differentiable almost everywhere,
that is, the set $M\setminus\textrm{Diff}(g)$ has measure zero.
If $M$ is infinite-dimensional, then the set of points of
differentiability of $g$, $\textrm{Diff}(g)$, is dense in $M$.
\end{thm}
\begin{proof}
Since $M$ is separable, it suffices to prove the result for any
small enough open set $U\subset M$ so that $g$ is Lipschitz on
$U$. Take a point $p\in U$. Since the exponential mapping at $p$
is locally almost an isometry, in particular Lipschitz, it
provides us with a chart $h=\Phi_{p}:V\to H$ which is Lipschitz,
for a suitably small open set $V\subset U$. Then the composition
$g\circ h^{-1}:h(V)\subset H\to\mathbb{R}$ is a Lipschitz function
from an open subset of a Hilbert space into $\mathbb{R}$.

In the case when $H$ is finite-dimensional, the classic theorem
of Rademacher tells us that $g\circ h^{-1}$ is differentiable
almost everywhere in $h(V)$ (see \cite{Federer}) and, since $h$
is a $C^1$ diffeomorphism (so $h$ preserves points of
differentiability and sets of measure zero), it follows that $g$
is differentiable almost everywhere in $V$.

If $H$ is infinite-dimensional then we can apply a celebrated
theorem of Preiss that ensures that every Lipschitz function from
an open set of an Asplund Banach space (such as the Hilbert space)
has at least one point of differentiability \cite{Preiss}. By
this theorem, it immediately follows that $g\circ h^{-1}$ is
differentiable on a dense subset of $h(V)$. Since again $h$ is a
$C^1$ diffeomorphism, we have that $g$ is differentiable on a
dense subset of $V$.

Finally, since $M$ can be covered by a countable union of such
open sets $V$ on each of which $g$ is Lipschitz, the result
follows.
\end{proof}
\begin{cor}
Let $M$ be a Riemannian manifold, and $f:M\to\mathbb{R}$ a convex
and locally bounded function. Then $f$ is differentiable on a
dense subset of $M$ (whose complement has measure zero if $M$ is
finite-dimensional).
\end{cor}
\begin{proof}
By Proposition \ref{every convex function is Lipschitz} we know
that $f$ is locally Lipschitz. Then, by Theorem \ref{Lipschitz
functions are differentiable somewhere}. it follows that $f$ is
differentiable on a dense subset of $M$.
\end{proof}

\medskip


\section{Hamilton-Jacobi equations in Riemannian manifolds}

\medskip

\noindent First order Hamilton-Jacobi equations are of the form
    $$
    F(x, u(x), du(x))=0
    $$
in the stationary case, and of the form
    $$
    F(t,x, u(x,t), du(t,x))=0
    $$
in the evolution case. These equations arise, for instance, in
optimal control theory, Lyapounov theory, and differential games.

Even in the simplest cases, such as the space $\mathbb{R}^{n}$, it
is well known that very natural Hamilton-Jacobi equations do not
always admit classical solutions. However, weaker solutions, such
as the so-called viscosity solutions, do exist under very general
assumptions. There is quite a large amount of literature about
viscosity solutions to Hamilton-Jacobi equations, see
\cite{Barles, C, CL1, CL2, CL3, CL4, CL5, CL6, CL7, CL8, CL9,
CL10, DGZ, DeGou} and the references cited therein, for instance.
All these works deal with Hamilton Jacobi equations in
$\mathbb{R}^{n}$ or in infinite-dimensional Banach spaces.

Examples of Hamilton-Jacobi equations also arise naturally in the
setting of Riemannian manifolds, see \cite{AbrahamMarsden}.
However, we do not know of any work that has studied nonsmooth
solutions, in general, or viscosity solutions, in particular, to
Hamilton-Jacobi equations defined on Riemannian manifolds (either
finite-dimensional or infinite-dimensional). This may be due to
the lack of a theory of nonsmooth calculus for functions defined
on Riemannian manifolds.

In this final section we will show how the subdifferential
calculus we have developed can be applied to get results on
existence and uniqueness of viscosity solutions to some
Hamilton-Jacobi equations defined on Riemannian manifolds. We
will also prove some results about ``regularity" (meaning
Lipschitzness) of viscosity solutions to some of these equations.

There are lots of Hamilton-Jacobi equations on Riemannian
manifolds $M$ for which the tools we have just developed could be
used in one way or the other to get interesting results about
viscosity solutions. For instance, one could get a maximum
principle for stationary first order Hamilton-Jacobi equations of
the type
    $$
    \left\{
\begin{matrix}u(x)+F(x, du(x))=0 \textrm{ for all } x\in\Omega \cr u(x)=0
\textrm{ for all } x\in \partial\Omega, \cr\end{matrix} \right.
    $$
where $\Omega$ is an open submanifold of $M$ with boundary
$\partial\Omega$. One could also prove a maximum principle for
parabolic Hamilton-Jacobi equations of the form
    $$
    \left\{
\begin{matrix}u_{t}+F(x, u_x)=0, \cr u(0,x)=u_{0}(x), \cr\end{matrix} \right.
    $$
where $u:\mathbb{R}^{+}\times M\to\mathbb{R}$, and
$u_{0}:M\to\mathbb{R}$ is an initial condition (assumed to be
bounded and uniformly continuous), in the manner of \cite[Section
6]{DeGou}.

However, this final section is only intended to give a glimpse of
the potential applications of nonsmooth calculus to the theory of
Hamilton-Jacobi equations on Riemannian manifolds, and not to
elaborate a comprehensive treatise on such equations. That is why
we will restrict ourselves mainly to one of the most interesting
examples of first-order Hamilton-Jacobi equations, namely
equations of the form $$(*) \left\{
\begin{matrix}u+G(du)=f \cr u \textrm{ bounded, } \cr\end{matrix} \right.$$
where $f:M\to \mathbb{R}$ is a bounded uniformly continuous
function, and $G:T^*M\to \mathbb{R}$ is a function defined on the
cotangent bundle of $M$. In fact these equations are really of the
form $$(*) \left\{
\begin{matrix}u+F(du)=0 \cr u \textrm{ bounded, } \cr\end{matrix} \right.$$
where $F:T^*M\to \mathbb{R}$, since we can always take a function
$F$ of the form $F(x,\xi_{x})=G(x,\xi_{x})-f(x)$.

\medskip

A bounded Fr{\'e}chet-differentiable function $u:M\to \mathbb{R}$
is a classical solution of the equation $(*)$ provided that
$$u(p)+F(p,du(p))=0 \ \ \ \hbox{ for every }\ p\in M.$$ Let us
now introduce the notion of viscosity solution.

\begin{defn}
{\em An upper semicontinuous (usc) function $u:M\to \mathbb{R}$\
is a {\em viscosity subsolution} of $u+F(du)=0$\ if
$u(p)+F(p,\zeta )\leq 0$\ for every $p\in M$\ and $\zeta \in
D^{+}u(p)$. A lower semicontinuous (lsc) function $u:M\to
\mathbb{R}$\ is a {\em viscosity supersolution} of $u+F(du)=0$\ if
$u(p)+F(p,\zeta )\geq 0$\ for every $p\in M$\ and $\zeta \in
D^{-}u(p)$. A continuous function $u:M\to \mathbb{R}$ is a {\em
viscosity solution} of $u+F(du)=0$ if it is both a viscosity
subsolution and a viscosity supersolution of $u+F(du)=0$.

We can define viscosity solutions on a open set $\Omega \subset
M$\ in a natural way when the functions are defined on $\Omega$.}
\end{defn}
\begin{rem}
{\em Since for a Fr{\'e}chet differentiable function $u$ we have
$D^{+}u(p)=D^{-}u(p)=\{ du(p) \}$, it is clear that every bounded
Frechet differentiable viscosity solution of $u+F(du)=0$ is a
classical solution of $(*)$.}
\end{rem}

\medskip

We are going to show the existence and uniqueness of viscosity
solutions to Hamilton-Jacobi equations of the form $(*)$ provided
that $F:T^*M\to \mathbb{R}$ is a function defined on the cotangent
bundle of $M$ which satisfies a certain uniform continuity
condition, see Definition \ref{intrinsically uniformly continuous
function} below. The manifold $M$ must also satisfy the following
requirement.

Throughout the remainder of this section $M$ will be a complete
Riemannian manifold (either finite-dimensional or
infinite-dimensional) such that $M$ satisfies conditions $(3)$ or
$(4)$ (which are both equivalent) of Proposition \ref{sufficient
conditions for M to be uniformly bumpable}, that is, $M$ is
uniformly locally convex and has a strictly positive injectivity
radius. Equivalently, there is a constant $r=r_{M}>0$ such that
for every $x\in M$ the mapping $\exp_{x}$ is defined on $B(0_{x},
r)\subset TM_{x}$ and provides a $C^\infty$ diffeomorphism
    $$
    \exp_{x}:B(0_{x}, r)\to B(x,r),
    $$
and the distance function is given by the expression
    $$
    d(y,x)=\|\exp_{x}^{-1}(y)\|_{x} \, \textrm{ for all }\, y\in
    B(x,r).
    $$
In particular, all compact manifolds satisfy this property. In
the remainder of this section the constant $r=r_{M}$ will be
fixed.

Note also that if $M$ satisfies condition $(3)$ of Proposition
\ref{sufficient conditions for M to be uniformly bumpable} then
$M$ is uniformly bumpable and therefore the Smooth Variational
Principle \ref{Smooth Variational Principle} holds for $M$.

We begin with a simple observation that if $M$ is uniformly
bumpable then so is $M\times M$.
\begin{lem}\label{the product is uniformly bumpable}
Let $M$ be a Riemannian manifold. If $M$ is uniformly bumpable
then $M\times M$ is uniformly bumpable as well.
\end{lem}
\begin{proof}
The natural Riemannian structure in $M\times M$ induced by $(M,
g)$ is the one given by
    $$
    (g\times g)_{(p_{1},p_{2})}\big( (v_{1},v_{2}),\, (w_{1},w_{2})\big):=
    g_{p_{1}}(v_{1}, w_{1})+g_{p_{2}}(v_{2}, w_{2}).
    $$
Let $d_{M\times M}$ denote the Riemannian distance that this
metric gives rise to in the product $M\times M$. It is obvious
that if $\gamma(t)=\big(\alpha(t), \beta(t)\big)$ is a path in
$M\times M$ then $\alpha$ and $\beta$ are paths in $M$ satisfying
    $$
    \max\{L(\alpha), L(\beta)\}\leq L(\gamma)\leq L(\alpha)+L(\beta)\leq
    2\max\{L(\alpha), L(\beta)\},
    $$
which implies that
    \begin{eqnarray*}
    & & \max\{d_{M}(x_{1}, y_{1}), d_{M}(x_{2}, y_{2})\}\leq
    d_{M\times M}\big( (x_{1}, x_{2}),\, (y_{1}, y_{2}) \big)\leq\\
    & & d_{M}(x_{1}, y_{1})+d_{M}(x_{2}, y_{2})
    \leq 2\max\{d_{M}(x_{1}, y_{1}), d_{M}(x_{2}, y_{2})\}
    \end{eqnarray*}
for every $x=(x_{1},x_{2}), y=(y_{1},y_{2})\in M\times M$.

Since $M$ is uniformly bumpable, there exist numbers $R=R_{M}>1$
and $r_{M}>0$ such that for every $p_{0}\in M$, $\delta\in (0,
r_{M})$ there exists a $C^1$ smooth function $b:M\to\mathbb
[0,1]$ such that $b(p_{0})=1$, $b(x)=0$ if
$d_{M}(x,p_{0})\geq\delta$, and $\sup_{x\in M}\|db(x)\|_{x}\leq
R/\delta$. Now take a point $p=(p_{1}, p_{2})\in M\times M$. For
any $\delta\in (0, r_{M})$, there are $C^1$ smooth bumps $b_{1},
b_{2}$ on $M$ such that $b_{i}(p_{i})=1$, $b_{i}(x_{i})=0$
whenever $d_{M}(x_{i},p_{i})\geq\delta$, and $\|db(x_{i})\|\leq
R/\delta$ for every $x_{i}\in M$; $i=1,2$. Define a $C^1$ smooth
bump $b:M\times M\to\mathbb{R}$ by
    $$
    b(x)=b(x_{1},x_{2})=b_{1}(x_{1}) b_{2}(x_{2}) \textrm{ for all }
    x=(x_{1}, x_{2})\in M\times M.
    $$
It is obvious that $b(p_{1}, p_{2})=1$. If $d_{M\times M}(x,p)\geq
2\delta$ we have that
    $$
    2\max\{d_{M}(x_{1}, p_{1}), d_{M}(x_{2}, p_{2})\}\geq
    d_{M\times M}(x,p)\geq 2\delta,
    $$
so $d_{M}(x_{i}, p_{i})\geq\delta$ for some $i\in\{1,2\}$, hence
$b_{i}(x_{i})=0$ for the same $i$, and $b(x)=0$. Finally, we have
that
    $$
    \| db(x_{1}, x_{2})\|_{(x_{1}, x_{2})}^{2}=
    \|db_{1}(x_{1})\|_{x_{1}}^{2}+\|db_{2}(x_{2})\|_{x_{2}}^{2}
    \leq 2(R/\delta)^{2},
    $$
which means that
    $$
    \|db(x)\|_{x}\leq \frac{2\sqrt{2} R}{2\delta}
    $$
for every $x=(x_{1}, x_{2})\in M\times M$. Therefore $M\times M$
satisfies the conditions in Definition \ref{uniformly bumpable
manifold}, with $R_{M\times M}=2\sqrt{2}R$, and $r_{M\times M}=2
r_{M}$.
\end{proof}

Since we are assuming that $M$ is uniformly locally convex and
$i(M)>r>0$, hence that the distance function $y\mapsto d(y,x)$ is
$C^\infty$ smooth on $B(x,r)\setminus\{x\}$ for every $x\in M$,
we can consider the distance function $d:M\times M\to\mathbb{R}$
and its partial derivatives $\partial d(x_{0}, y_{0})/\partial x$
and $\partial d(x_{0}, y_{0})/\partial y$. We next see that these
partial derivatives satisfy a nice antisymmetry property. In
order to compare them in a natural way we need to use the
parallel translation from $TM_{x_{0}}$ to $TM_{y_{0}}$ along the
geodesic joining $x_{0}$ to $y_{0}$ (note that there is a unique
minimizing geodesic joining $x_{0}$ to $y_{0}$ because $M$ is
uniformly locally convex and $d(x_{0}, y_{0})<r$).
\begin{notation}
{\em Let $x_{0}, y_{0}\in M$ be such that $d(x_{0}, y_{0})<r$. Let
$\gamma(t)=\exp_{x_{0}}(tv_{0})$, $0\leq t\leq 1$ be the unique
minimizing geodesic joining these two points. For every vector
$w\in TM_{x_{0}}$, we denote
    $$
    L_{x_{0}y_{0}}(w)=P_{0,\gamma}^{1}(w)
    $$
the parallel translation of $w$ from $x_{0}$ to $y_{0}$ along
$\gamma$. Recall that the mapping $L_{x_{0}y_{0}}:TM_{x_{0}}\to
TM_{y_{0}}$ is a linear isometry, with inverse
$L_{y_{0}x_{0}}:TM_{y_{0}}\to TM_{x_{0}}$. As we customarily
identify $TM_{p}$ with $T^{*}M_{p}$ (via the linear isometry
$v\mapsto \langle v, \cdot\rangle_{p}$), the isometry
$L_{x_{0}y_{0}}$ induces another linear isometry between the
cotangent fibers $T^{*}M_{x_{0}}$ and $T^{*}M_{y_{0}}$. We will
still denote this new isometry by
$L_{x_{0}y_{0}}:T^{*}M_{x_{0}}\to T^{*}M_{y_{0}}$.}
\end{notation}
\begin{lem}\label{antisymmetry of the partial derivatives of the distance}
Let $x_{0}, y_{0}\in M$ be such that $0<d(x_{0}, y_{0})<r$. Then
    $$
    L_{y_{0}x_{0}}\Big(\frac{\partial d(y_{0}, x_{0})}{\partial
    y}\Big)=\,
    -\frac{\partial d(x_{0}, y_{0})}{\partial x}.
    $$
\end{lem}
\begin{proof}
Denote $r_{0}=d(x_{0},y_{0})<r$. Consider the geodesic
$\gamma(t)=\exp_{x_{0}}(tv_{0})$, $0\leq t\leq 1$, where
$y_{0}=\exp_{x_{0}}(v_{0})$. By the definitions of parallel
translation and geodesic it is clear that
    $$
    L_{x_{0}y_{0}}(v_{0})=\gamma'(1)=d\exp_{x_{0}}(v_{0})(v_{0}).
    $$
On the other hand, under the current assumptions on $M$, and by
the Gauss Lemma (see \cite{Klingenberg, Lang}), we know that
$\gamma'(1)$ is orthogonal to the sphere $S(x_{0}, r_{0})=\{y\in
M : d(y, x_{0})=r_{0}\}=\exp_{x_{0}}(S(0_{x_{0}}, r_{0}))$. Since
this sphere $S(x_{0}, r_{0})$ is a one-codimensional submanifold
of $M$ defined as the set of zeros of the smooth function
$y\mapsto d(y, x_{0})-r_{0}$ and (as is easily checked)
$$\frac{\partial d(y_{0},x_{0})}{\partial y}\neq 0,$$ we also
have that this partial derivative is orthogonal to the sphere
$S(x_{0}, r_{0})$ at the point $y_{0}$. Therefore,
    $$
    L_{x_{0}y_{0}}(v_{0})=
    \gamma'(1)=
    \lambda\frac{\partial d(y_{0},x_{0})}{\partial y}
    $$
for some $\lambda\neq 0$. Furthermore, since the function
$t\mapsto d(\gamma(t), x_{0})$ is increasing, we get that
$\lambda>0$. Finally, it is clear that $y\mapsto d(y,x_{0})$ is
$1$-Lipschitz, and
    $$
    \|\frac{\partial d(y_{0},x_{0})}{\partial y}\|=1,
    $$
from which we deduce that
$\lambda=\|L_{x_{0}y_{0}}(v_{0})\|_{y_{0}}=\|v_{0}\|_{x_{0}}$, and
    $$
    L_{x_{0}y_{0}}(v_{0})=
    \|v_{0}\|_{x_{0}}
    \frac{\partial d(x_{0},y_{0})}{\partial y}.
    \eqno(1)
    $$
Now consider the geodesic $\beta(t)=\exp_{y_{0}}(t w_{0})$,
$0\leq t\leq 1$, where $\exp_{y_{0}}(w_{0})=x_{0}$. By the
definitions of parallel translation and geodesic we know that
    $$
    L_{x_{0}y_{0}}(v_{0})=-w_{0},\, \textrm{ and } \,\,
    \|w_{0}\|_{y_{0}}=\|v_{0}\|_{x_{0}}.
    \eqno(2)
    $$
A completely analogous argument to the one we used for $\gamma$
above shows that
    $$
    L_{y_{0}x_{0}}(w_{0})=
    \|w_{0}\|
    \frac{\partial d(x_{0},y_{0})}{\partial x}. \eqno(3)
    $$
By combining $(1)$, $(2)$ and $(3)$ we immediately get that
    $$
    L_{y_{0}x_{0}}\Big(\frac{\partial d(y_{0}, x_{0})}{\partial
    y}\Big)=\frac{v_{0}}{\|v_{0}\|_{x_{0}}}=
    -\frac{L_{y_{0}x_{0}}(w_{0})}{\|w_{0}\|_{y_{0}}}=
    -\frac{\partial d(x_{0}, y_{0})}{\partial x}.
    $$
\end{proof}
The following Proposition can be viewed as a perturbed
minimization principle for the sum or the difference of two
functions. Its proof is a consequence of the Smooth Variational
Principle \ref{Smooth Variational Principle} and Lemma
\ref{antisymmetry of the partial derivatives of the distance}
\begin{prop}\label{HC prop 1}
Let $u,v:M\to \mathbb{R}$ be two bounded functions which are
upper semicontinuous (usc) and lower semicontinuous (lsc)
respectively. Then, for every $\varepsilon >0$, there exist
$x_{0}, y_{0}\in M$, and $\zeta \in D^{+}u(x_{0})$, $\xi \in
D^{-}v(y_{0})$ such that
\begin{enumerate}
\item[{(i)}] $d(x_{0},y_{0})<\varepsilon$,
\item[{(ii)}] $\|\zeta -L_{y_{0}x_{0}}(\xi)\|_{x_{0}}<\varepsilon$
\item[{(iii)}] $v(z)-u(z)\geq v(y_{0})-u(x_{0})-\varepsilon$ for every $z\in M$.
\end{enumerate}
{\em Here $L_{y_{0}x_{0}}:T^{*}M_{y_{0}}\to T^{*}M_{x_{0}}$ stands
for the parallel translation.}
\end{prop}
\begin{proof}
We can obviously assume that $\varepsilon<r(M)$. Let
$b:\mathbb{R}\to\mathbb{R}$ be a $C^\infty$ smooth function such
that $b$ is non-increasing,
    $$
    b(t)=b(0)>2(\|v\|_{\infty}+\|u\|_{\infty})+\varepsilon \, \textrm{
    if } t\leq \varepsilon/4, \, \textrm{ and }\, b(t)=0\,
    \textrm{ if }\, t\geq\varepsilon. \eqno(1)
    $$
Define the function $w: M\times M\to\mathbb{R}$ by
    $$
    w(x,y)=v(y)-u(x)-b(d(x,y))\, \textrm{ for all }\, (x,y)\in M\times M.
    $$
The function $w$ is lower semicontinuous and bounded. By Lemma
\ref{the product is uniformly bumpable} we know that $M\times M$
is uniformly bumpable, and $M\times M$ is obviously complete, so
we can apply the Smooth Variational Principle \ref{Smooth
Variational Principle} to the function $w$ to get a pair $(x_{0},
y_{0})\in M\times M$ and a $C^1$ smooth function $g:M\times
M\to\mathbb{R}$ such that
\begin{enumerate}
\item[{(a)}] $\|g\|_{\infty}<\varepsilon/2 >\|dg\|_{\infty}$
\item[{(b)}] $v(y)-u(x)-b(d(x,y))-g(x,y)\geq
v(y_{0})-u(x_{0})-b(d(x_{0}, y_{0})) -g(x_{0}, y_{0})$ for all
$x,y\in M$.
\end{enumerate}
If we take $x=x_{0}$ in $(b)$ we get that $v$ is
subdifferentiable at the point $y_{0}$, and
    $$
    \xi:=\frac{\partial g(x_{0},y_{0})}{\partial y}
    +\frac{\partial (b\circ d)(x_{0},y_{0})}{\partial y}
    \in D^{-}v(y_{0}). \eqno(2)
    $$
In a similar manner, by taking $y=y_{0}$ in $(b)$ we get that
    $$
    \zeta:=-\Big(\frac{\partial g(x_{0},y_{0})}{\partial x}
    +\frac{\partial (b\circ d)(x_{0},y_{0})}{\partial x}\Big)
    \in D^{+}u(x_{0}). \eqno(3)
    $$
Let us note that
\begin{eqnarray*}
& & L_{y_{0}x_{0}}\Big(\frac{\partial (b\circ
d)(x_{0},y_{0})}{\partial y}\Big)+\frac{\partial (b\circ
d)(x_{0},y_{0}))}{\partial x}=\\ & & b'(d(x_{0},
y_{0}))\Big[L_{y_{0}x_{0}}\Big(\frac{\partial
d(x_{0},y_{0})}{\partial y}\Big)+\frac{\partial
d(x_{0},y_{0})}{\partial x}\Big]=0, \hspace{3.7cm} (4)
\end{eqnarray*}
thanks to Lemma \ref{antisymmetry of the partial derivatives of
the distance} when $x_{0}\neq y_{0}$, and to the definition of
$b$ when $x_{0}=y_{0}$. Therefore,
\begin{eqnarray*}
& &\|L_{y_{0}x_{0}}(\xi)-\zeta\|_{x_{0}}=\\ & &
\|L_{y_{0}x_{0}}\Big(\frac{\partial g(x_{0},y_{0})}{\partial
y}\Big) +L_{y_{0}x_{0}}\Big(\frac{\partial (b\circ
d)(x_{0},y_{0})}{\partial y}\Big) +\frac{\partial
g(x_{0},y_{0})}{\partial x}
    +\frac{\partial (b\circ d)(x_{0},y_{0})}{\partial x}\|_{x_{0}}=\\
& &\|L_{y_{0}x_{0}}\Big(\frac{\partial g(x_{0},y_{0})}{\partial
y}\Big) +\frac{\partial g(x_{0},y_{0})}{\partial
x}\|_{x_{0}}\leq\\ & &\|\frac{\partial g(x_{0},y_{0})}{\partial
y}\|_{y_{0}} +\|\frac{\partial g(x_{0},y_{0})}{\partial
x}\|_{x_{0}}\leq
\|dg\|_{\infty}+\|dg\|_{\infty}<\frac{\varepsilon}{2}+
\frac{\varepsilon}{2}=\varepsilon,
\end{eqnarray*}
which shows $(ii)$.

On the other hand, if we had $d(x_{0}, y_{0})\geq\varepsilon$
then, by taking $x=y=z$ in $(b)$ we would get
\begin{eqnarray*}
& & b(0)\leq
v(z)-u(z)-g(z,z)+g(x_{0},y_{0})-v(y_{0})+u(x_{0})\leq\\ & &
2(\|v\|_{\infty}+\|u\|_{\infty})+\varepsilon,
\end{eqnarray*}
which contradicts the definition of $b$, see $(1)$ above.
Therefore $d(x_{0},y_{0})<\varepsilon$ and $(i)$ is proved.

Finally, if we take $z=x=y$ in $(b)$ and we bear in mind that
$\|g\|_{\infty}<\varepsilon/2$ and the function $b$ is
non-increasing, we get that
\begin{eqnarray*}
& & v(z)-u(z)\geq v(y_{0})-u(x_{0})+b(0)-b(d(x_{0},
y_{0}))+g(z,z)-g(x_{0}, y_{0})\geq\\ & &
v(y_{0})-u(x_{0})+0-\varepsilon/2-\varepsilon/2=
v(y_{0})-u(x_{0})-\varepsilon,
\end{eqnarray*}
which shows $(iii)$ and finishes the proof.
\end{proof}
\begin{rem}
{\em Let us observe that the preceding proposition is no longer
true if the manifold is not complete. For example:
$M=(0,1)\subset\mathbb{R}$, $g(x)=x$, $f(x)=0$, and
$\varepsilon>0$ small.}
\end{rem}
\begin{defn}
{\em For a given open set $\Omega \subset M$\ and a function
$u:\Omega \to \mathbb{R}$, we define the upper semicontinuous
envelope of $u$, which we denote $u^{*}$, by
    $$
    u^{*}(x)=\inf\{ v(x) \, |\, v:\Omega\to\mathbb{R}
    \textrm{ is continuous and } u\leq v \textrm{ on } \Omega\}
    \,\, \textrm{ for any } x\in \Omega.
    $$
In a similar way we define the lower semicontinuous envelope,
denoted by $u_{*}$.}
\end{defn}

\begin{prop}\label{prop 69}
Let $\Omega$ be an open subset of $M$. Let $\mathcal F$\ be  a
uniformly bounded family of upper semicontinuous functions on
$\Omega$, and let $u=\sup\{ v:v\in \mathcal F\}$. Then, for every
$p\in \Omega$, and every $\zeta \in D^+u^*(p)$, there exist
sequences $\{ v_n\} \subset \mathcal F$, and $\{ (p_n, \zeta _n)\}
\subset T^*(\Omega )$, with $\zeta _n\in D^+v_n(p_n)$ and such
that
\begin{enumerate}
\item $\lim_nv_n(p_n)=u^*(p)$
\item $\lim_n(p_n,\zeta _n)=(p,\zeta )$.
\end{enumerate}
\end{prop}
\begin{proof}
Fix a chart $(U,\varphi )$, with $p\in U$. Let us consider the
family $\mathcal F \circ \varphi ^{-1}=\{ v\circ \varphi
^{-1}:v\in \mathcal F\}$. The functions of this collection are
upper semicontinuous on $\varphi (U\cap \Omega)$, and the family
is uniformly bounded. On the other hand $u\circ \varphi
^{-1}=sup\{ v\circ \varphi ^{-1}:v\in \mathcal F\}$, and
$u^*\circ \varphi ^{-1}=(u\circ \varphi ^{-1})^*$.  Now apply
\cite[Proposition VIII.1.6]{DGZ} (which is nothing but the result
we want to prove in the case of a Banach space) to the Hilbert
space, the open set $\varphi (U\cap \Omega)$, the family
$\mathcal{F} \circ \varphi ^{-1}$, the point $\varphi (p)$, and
the superdifferential $\zeta \circ d\varphi (p)^{-1}\in D^+(u\circ
\varphi ^{-1})(\varphi (p))$. We get sequences $\{ \varphi
(p_n)\}$\ in $\varphi (U\cap \Omega)$, $\{ v_n\circ \varphi
^{-1}\}$\ in $\mathcal{F} \circ \varphi ^{-1}$, and $\zeta _n\circ
d\varphi (p)^{-1}\in D^+(v_n\circ \varphi ^{-1})(\varphi (p_n))$
such that $\lim_n\varphi (p_n)=\varphi (p)$, $\lim_n\zeta _n\circ
d\varphi (p)^{-1}=\zeta \circ d\varphi (p)^{-1}$, and
$$\lim_n(v_n\circ \varphi ^{-1})(\varphi (p_n))=(u\circ \varphi
^{-1})^*(\varphi (p))=u^*\circ \varphi ^{-1}(\varphi (p)).$$ Hence
$\lim_np_n=p$, $\lim_nv_n(p_n)=u^*(p)$, and $$\lim_n\zeta _n\circ
d\varphi (p_n)^{-1}=\lim_n\zeta _n\circ d\varphi (p_n)^{-1}\circ
d\varphi (p) \circ d\varphi (p)^{-1}=\lim_n\zeta _n\circ d\varphi
(p)^{-1}$$ because $\varphi$ is $C^1$, so $\lim_nd\varphi
(p_n)^{-1}\circ d\varphi (p)=\textrm{id}$. The result follows
trivially from the local representation of the cotangent bundle.
\end{proof}

Now we introduce the notion of uniform continuity that we have to
require of $F:T^{*}M\to\mathbb{R}$ in order to prove the existence
and uniqueness of viscosity solutions to the Hamilton-Jacobi
equation $(*)$.

\begin{defn}\label{intrinsically uniformly continuous function}
{\em We will say that a function $F:T^{*}M\to\mathbb{R}$ is
intrinsically uniformly continuous provided that for every
$\varepsilon>0$ there exists $\delta\in (0, r_{M})$ such that
    $$
    d(x,y)\leq\delta, \, \zeta\in T^{*}M_{x}, \, \xi\in
    T^{*}M_{y}, \, \|\zeta-L_{yx}(\xi)\|_{x}\leq\delta
    \implies |F(x,\zeta)-F(y,\xi)|\leq\varepsilon.
    $$
    }
\end{defn}
\begin{rem}
{\em It should be noted that if $F$ satisfies the above definition
then $F$ is continuous. This is obvious once we notice that the
mapping
    $$
    (x,\zeta)\in T^{*}M_{x}\mapsto L_{x x_{0}}(\zeta)
    $$
is continuous at $(x_{0},\zeta_{0})$, that is, if
$(x_{n},\zeta_{n})\to (x_{0},\zeta_{0})$ in $T^{*}M$ then
$L_{x_{n} x_{0}}(\zeta_{n})\to\zeta_{0}$, \, for every
$(x_{0},\zeta_{0})\in T^{*}M$. The fact that this mapping is
continuous is an easy consequence of the definition of the
parallel translation along a curve as a solution to an ordinary
linear differential equation.}
\end{rem}
\begin{rem}\label{usual uniform continuity implies intrinsic
uniform continuity} {\em Consider a finite-dimensional manifold
$M$ embedded in $\mathbb{R}^{n}$, so $T^{*}M\subset
\mathbb{R}^{2n}$. Assume that $M$ satisfies the following
condition:
\begin{center}
$\forall\varepsilon$\, $\exists$\, $\delta>0$: $v, h\in TM_{x}$\,
$x\in M$, $\|v\|_{x}\leq\delta$\, $\implies$
$\|d\exp_{x}(v)(h)-h\|_{x}\leq\varepsilon$
\end{center}
(note in particular that this condition is automatically met when
$M$ is compact, and in many other natural examples). Then every
function $F:T^{*}M\to\mathbb{R}$ which is uniformly continuous
with respect to the usual Euclidean metric in $\mathbb{R}^{2n}$ is
intrinsically uniformly continuous as well, as is easily seen.
Consequently there are lots of natural examples of intrinsically
uniformly continuous functions $F:T^{*}M\to\mathbb{R}$.}
\end{rem}

Now we can prove the following {\em maximum principle} for
Hamilton-Jacobi equations of the form $(*)$.
\begin{thm}\label{maximum principle}
Let $f,g:M\to \mathbb{R}$\  be bounded uniformly continuous
functions, and $F:T^*M\to \mathbb{R}$ be intrinsically uniformly
continuous. If $u$\ is a bounded viscosity subsolution of
$u+F(du)=f$\ and $v$\ is a bounded viscosity supersolution of
$v+F(dv)=g$, then $v-u\geq \inf (g-f)$.
\end{thm}
\begin{proof}
If $\varepsilon >0$\ is given, then, by Proposition \ref{HC prop
1}, there exist $p, q\in U$, and $\zeta \in D^+u(p)$, $\xi \in
D^-v(q)$ such that
\item{i)} $d(p,q)<\varepsilon$, $\|\zeta -L_{qp}(\xi)\|_{p}<\varepsilon$
\item{ii)} $v(x)-u(x)\geq v(q)-u(p)-\varepsilon$\ for every $x\in M$.

Since $u$\ and $v$\ are viscosity sub and super solutions
respectively, we have $u(p)+F(p,\zeta )\leq f(p)$\ and
$v(q)+F(q,\xi )\geq g(q)$. Hence, for every $x\in M$,
$$v(x)-u(x)\geq v(q)-u(p)-\varepsilon \geq g(q)-F(q,\xi
)-f(p)+F(p,\zeta )-\varepsilon \geq$$ $$\geq  \inf
(g-f)+(f(q)-f(p))+(F(p,\zeta )-F(q,\xi ))-\varepsilon. $$ Now, if
we let $\varepsilon \to 0^+$, we have that $f(q)-f(p)$ goes to $0$
because $f$\ is uniformly continuous. On the other hand, the fact
that $F$ is intrinsically uniformly continuous implies that
$F(p,\zeta )-F(q,\xi )$ goes to $0$ as $\varepsilon\to 0^{+}$.
Consequently we obtain $v-u\geq \inf (g-f)$.
\end{proof}

\begin{prop}\label{prop 613}
Let $\Omega$ be an open subset of $M$. Let $\mathcal F$ be  a
uniformly bounded family of functions on $\Omega$, and let
$u=\sup\{ v:v\in \mathcal F\}$. If every $v$ is a viscosity
subsolution of $u+F(du)=0$, then $u^{*}$ is a viscosity
subsolution of $u+F(du)=0$.
\end{prop}
\begin{proof}
Let $p\in \Omega$\ and $\zeta \in D^+u^*(p)$. According to
Proposition \ref{prop 69}, there exist sequences $\{ v_n\} \subset
\mathcal F$, and $\{ (p_n, \zeta _n)\} \subset T^*(\Omega)$ with
$\zeta _n\in D^{+}v_n(p_n)$ and such that
\item{i)} $\lim_nv_n(p_n)=u^{*}(p)$
\item{ii)} $\lim_n(p_n,\zeta _n)=(p,\zeta)$.

Since $v_n$\ are viscosity subsolutions of $u+F(du)=0$, we have
that $v_n(p_n)+F(p_n,\zeta _n)\leq 0$. Hence $u^{*}(p)+F(p,\zeta
)\leq 0$.
\end{proof}

\begin{cor}
The supremum of two viscosity subsolutions is a viscosity
subsolution.
\end{cor}

\begin{thm}\label{existence and uniqueness of viscosity solutions to HJ}
Let $M$ be a complete Riemannian manifold which is uniformly
locally convex and has a strictly positive injectivity radius. Let
$F:T^*M\to \mathbb{R}$ be an intrinsically uniformly continuous
function (see Definition \ref{intrinsically uniformly continuous
function} above). Assume also that there is a constant $A>0$ so
that $-A\leq F(x,0_{x})\leq A$\ for every $x\in M$. Then, there
exists a unique bounded viscosity solution of the equation
$u+F(du)=0$.
\end{thm}
\begin{proof}
Uniqueness follows from Theorem \ref{maximum principle}, by taking
$f=g=0$. In order to show existence, let us define $\mathcal{F}$
as the family of the viscosity subsolutions $w:M\to \mathbb{R}$ to
$u+F(du)=0$ that satisfy $$-A\leq w(p) \leq A\ \ \ \hbox{ for
every }\ p\in M.$$ The family $\mathcal F$ is nonempty, as the
function $w_{0}(p)=-A$ belongs to $\mathcal{F}$ (because
$-A+F(p,0_{p})\leq 0$). Let $u$ be the upper semicontinuous
envelope of $\sup \{ w: w\in \mathcal{F}\}$, and $v$ be the lower
semicontinuous envelope of $u$. By the definition, we have $v\leq
u$. On the other hand, according to Proposition \ref{prop 613},
$u$ is a viscosity subsolution of $u+F(du)=0$.
\begin{claim}$
v$ is a viscosity supersolution of $u+F(du)=0$.
\end{claim}
\noindent Once we have proved the claim, we have that $u\leq v$ by
Proposition \ref{maximum principle}, hence $u=v$ is a viscosity
solution, and existence is established.

So let us prove the claim. If $v$\ is not a viscosity
supersolution, there exist $p_0\in M$\ and $\zeta _0\in D^-v(p)$\
such that $v(p_0)+F(p_0,\zeta _0)<0$. By Theorem \ref{equivalent
definitions of subdifferential}(5), there exists a $C^1$ smooth
function $h:M\to \mathbb{R}$ with $\zeta_{0}=dh(p_0)$ and such
that $v-h$ attains a global minimum at $p_0$. Hence we may assume
that $$v(p_0)+F(p_0,dh(p_0))<0,\ \ v(p_0)=h(p_0), \hbox{ and }
h(p)\leq v(p) \hbox{ for all } p\in M. \eqno(1)$$ From the
inequality $h(p)\leq v(p)\leq u(p)\leq A$ we get $h(p_0)<A$:
otherwise $A-h$ would have a local minimum at $p_0$, and
consequently $dh(p_0)=0$, which implies
$v(p_0)+F(p_0,dh(p_0))=h(p_0)+F(p_0,dh(p_0))=A+F(p_{0},0)\geq A-A
=0$, a contradiction with $(1)$.

Now we can take a number $\delta >0$ and a $C^1$ smooth function
$b:M\to [0,\infty)$ with support on $B(p_0,\delta )$, $b(p_0)>0$,
and such that $||b||_{\infty}$, $||db||_{\infty}$ are small enough
so that $$h(p)+b(p)+F(p,dh(p)+db(p))<0 \ \ \hbox{ for every } \ \
p\in B(p_{0},2\delta), \textrm{ and } \eqno(2)$$ $$h(p)+b(p)\leq
A \ \ \hbox{ for every } \ \ p\in M \eqno(3).$$ This is possible
because of $(1)$ and the fact that $F$ is continuous.

Let us consider the following function: $$w(p)=\left\{
\begin{matrix}
\max\{h(p)+b(p),u(p)\} \ \ \ \ \ \hbox{ if }\ \ \ p\in
B(p_0,2\delta ) \cr u(p)\ \ \ \ \ \ \ \ \ \ \ \ \ \ \ \ \ \ \ \
\hbox{ otherwise. } \cr
\end{matrix}
\right.$$ We have that $w(p)=u(p)$ if $p\in \Omega_{1}:=M\setminus
\overline B (p_0,\delta )$, because $u(p)\geq v(p)\geq
h(p)=h(p)+b(p)$ whenever $p\in B (p_0,2\delta )\setminus \overline
B (p_0,\delta )$. Therefore $w$ is a viscosity subsolution of
$u+F(du)=0$ on $\Omega_{1}$. On the other hand, bearing $(2)$ in
mind, it is clear that $w$ is the maximum of two viscosity
subsolutions on $\Omega_{2}:=B (p_0,2\delta)$, and consequently
$w$ is a viscosity subsolution on $\Omega_{2}$. Therefore $w$ is a
viscosity subsolution of $u+F(du)=0$ on $M=\Omega_{1}\cup
\Omega_{2}$. This implies that $w\in \mathcal{F}$, since $-A\leq
u\leq w$ and $w(p)\leq A$, by $(3)$.

Finally, we have that $u\geq w$, because $u\geq sup \mathcal{F}$.
Therefore we have $u(p)\geq w(p)\geq h(p)+b(p)$\ on $B(p_0,\delta
)$, and in particular $v(p_0)=u_*(p_0)\geq h(p_0)+b(p_0)>h(p_0)$,
which contradicts $(1)$.
\end{proof}
When $M$ is compact, the preceding Theorem takes on a simpler
appearance.
\begin{cor}\label{HJ exuni for compact manifolds}
Let $M$ be a compact Riemannian manifold, $f:M\to\mathbb{R}$ a
continuous function, and $F:T^{*}M\to \mathbb{R}$ an intrinsically
uniformly continuous function. Then there exists a unique
viscosity solution of the equation $u+F(du)=f$.
\end{cor}
\begin{proof}
This follows immediately from Theorem \ref{existence and
uniqueness of viscosity solutions to HJ}, taking into account the
following facts: 1) if $M$ is compact then $M$ is uniformly
locally convex and $i(M)>0$ (see Remarks \ref{continuity of the
convexity radius} and \ref{continuity of the injectivity radius}
above); 2) every viscosity solution $u$ is continuous, hence $u$
is bounded on the compact manifold $M$; and of course 3) $f$ is
uniformly continuous because $f$ is continuous on $M$, compact.
\end{proof}
\begin{rem}
{\em In particular, when a compact manifold $M$ is regarded as
embedded in $\mathbb{R}^{n}$, so $T^{*}M\subset \mathbb{R}^{2n}$,
and $F:T^{*}M\to\mathbb{R}$ is uniformly continuous with respect
to the usual Euclidean metric in $\mathbb{R}^{2n}$, then Corollary
\ref{HJ exuni for compact manifolds} and Remark \ref{usual uniform
continuity implies intrinsic uniform continuity} yield the
existence of a unique viscosity solution to the equation
$u+F(du)=f$.}
\end{rem}
However, the requirement that $F$ is uniformly continuous cannot
be relaxed in principle, because the cotangent bundle $T^{*}M$ is
never compact, so, even though $F$ is continuous, we cannot ensure
that $F$ is uniformly continuous on $T^{*}M$.

\begin{rem}\label{how to dispense with the conditions on M}
{\em It should be noted that one may pose a Hamilton-Jacobi
problem such as $$u+F(du)=0, \,\, u \textrm{ bounded } \eqno(*)$$
on a manifold $M$ without presupposing any Riemannian structure
defined on $M$. Then one may consider the natural question whether
it is possible to find a suitable Riemannian structure $g$ which
makes $M$ uniformly locally convex and with a positive injectivity
radius, and which makes $F$ intrinsically uniformly continuous as
well. In other words, one can seek for a Riemannian manifold $N$
with positive convexity and injectivity radii and a diffeomorphism
$\psi:N\to M$ so that the function $G:T^{*}N\to\mathbb{R}$ is
intrinsically uniformly continuous, where $G=F\circ (T^{*}\psi)$,
with $T^{*}\psi:T^{*}N\to T^{*}M$ defined by
$T^{*}\psi(x,\eta):=\big(\psi(x),\eta\circ
(d\psi^{-1})(\psi(x))\big)$. Then, by Theorem \ref{existence and
uniqueness of viscosity solutions to HJ}, the equation
    $$
    v+G(dv)=0, \,\, v \textrm{ bounded } \eqno(**)$$
has a unique viscosity solution. But it is obvious that $v$ is a
viscosity solution to $(**)$ if and only if the function
$u=v\circ\psi^{-1}$ is a viscosity solution to $(*)$. Hence the
equation $(*)$ has a unique viscosity solution as well. This means
that Theorem \ref{existence and uniqueness of viscosity solutions
to HJ} above is applicable to even more situations than one might
think of at a single glance. The following example reveals the
power of this scheme.}
\end{rem}
\begin{ex}
{\em Let $M$ be the submanifold of $\mathbb{R}^{3}$ defined by
$$z=\frac{1}{x^{2}+y^{2}},$$ let
$F:T^{*}M\subset\mathbb{R}^{6}\to\mathbb{R}$, and consider the
Hamilton-Jacobi equation $u+F(du)=0$. If we endow $M$ with its
natural Riemannian structure inherited from $\mathbb{R}^{3}$, $M$
will not be uniformly locally convex, and besides $i(M)=0$, so
Theorem \ref{existence and uniqueness of viscosity solutions to
HJ} is not directly applicable. Now let us define $N$ by
    $$
    z=\frac{1}{x^{2}+y^{2}-1}, \,\,\, z>0,
    $$
with its natural Riemannian structure as a submanifold of
$\mathbb{R}^{3}$. It is clear that $N$ is uniformly locally convex
and has a positive injectivity radius. The mapping $\psi:N\to M$
defined by
    $$
    F(x,y,z)=\Big(\frac{\sqrt{x^{2}+y^{2}-1}}{\sqrt{x^{2}+y^{2}}}\,
    x,\, \frac{\sqrt{x^{2}+y^{2}-1}}{\sqrt{x^{2}+y^{2}}}\, y,\,
    z\Big)
     $$
is a $C^\infty$ diffeomorphism. Assume that the function $G=F\circ
(T^{*}\psi):T^{*}N\to\mathbb{R}$ is uniformly continuous with
respect to the usual metric in $\mathbb{R}^{6}$. Since $N$
satisfies the property of Remark \ref{usual uniform continuity
implies intrinsic uniform continuity}, we have that $G$ is
intrinsically uniformly continuous. Therefore, by the preceding
Remark \ref{how to dispense with the conditions on M}, the
equation $u+F(du)=0$, $u$ bounded on $M$, has a unique viscosity
solution.}
\end{ex}

\medskip

Now let us see how Deville's mean value Theorem \ref{TVM Deville}
allows to deduce a result on the regularity of viscosity
solutions (or even subsolutions) to Hamilton-Jacobi equations
with a ``coercive" structure.

\begin{cor}\label{regularity of viscosity solutions}
Let $M$ be a Riemannian manifold, and $F:\mathbb{R}\times
T^{*}M\to\mathbb{R}$ a function. Consider the following
Hamilton-Jacobi equation:
    $$
    F(u(x), du(x))=0. \eqno(\textrm{{\em HJ3}})
    $$
Assume that there exists a constant $K>0$ such that $F(t,
\zeta_{x})>0$ whenever $\|\zeta_{x}\|_{x}\geq K$ and
$t\in\mathbb{R}$. Let $u:M\to\mathbb{R}$ be a viscosity
subsolution of {\em (HJ3)}. Then:
\begin{enumerate}
\item $u$ is $K$-Lipschitz, that is, $|u(x)-u(y)|\leq K d(x,y)$
for every $x,y\in M$.
\item If $M$ is finite dimensional, $u$ is Fr{\'e}chet
differentiable almost everywhere.
\item If $M$ is infinite-dimensional $u$ is Fr{\'e}chet
differentiable on a dense subset of $M$.
\end{enumerate}
\end{cor}
\begin{proof}
If $u$ is a viscosity subsolution then $F(u(x),\zeta_{x})\leq 0$
for every $x\in M$ and $\zeta_{x}\in D^{-}u(x)$. Hence
$\|\zeta_{x}\|\leq K$ for every  $\zeta_{x}\in D^{-}u(x)$
(otherwise $F(u(x),\zeta_{x})>0$, a contradiction). Then, by
Theorem \ref{TVM Deville}, $u$ is $K$-Lipschitz.

On the other hand, $(2)$ and $(3)$ follow immediately from
Theorem \ref{Lipschitz functions are differentiable somewhere}.
\end{proof}

\medskip

Let us conclude with a brief study of a HJ equation which is not
of the form $(*)$ above, but which is still very interesting
because of the geometrical significance of its unique viscosity
solution. Let $M$ be a complete Riemannian manifold, $\Omega$ a
bounded open subset of $M$, and let $\partial\Omega$ be the
boundary of $\Omega$. Consider the Hamilton-Jacobi equation
    $$
    (\textrm{HJ4}) \left\{\begin{matrix} \|du(x)\|_{x}=1 \textrm{ for all }
    x\in\Omega \cr u(x)=0 \textrm{ for all } x\in \partial\Omega.
    \cr\end{matrix} \right.
    $$
There is no classical solution of (HJ4). Indeed, if we had a
function $u:\overline{\Omega}\subset M\to\mathbb{R}$ which is
differentiable on $\Omega$ and satisfies $\|du(x)\|_{x}=1$ for
$x\in\Omega$ and $u=0$ on $\partial\Omega$, then we could apply
Theorem \ref{Rolle for manifolds} to find a point $x_{0}\in\Omega$
so that $\|du(x_{0})\|_{x_{0}}<1/2$, a contradiction.

Nevertheless, we are going to see that there is a unique
viscosity solution to (HJ4), namely the distance function to the
boundary $\partial\Omega$. By definition, a function $u$ is a
viscosity solution to (HJ4) if and only if $u$ is continuous;
$u=0$ on $\partial\Omega$; $\|\zeta\|_{x}\geq 1$ for all $\zeta\in
D^{-}u(x)$, $x\in\Omega$; and $\|\zeta\|_{x}\leq 1$ for all
$\zeta\in D^{+}u(x)$, $x\in\Omega$.

\begin{thm}
Let $M$ be a complete Riemannian manifold, and $\Omega$ a bounded
open subset of $M$ with boundary $\partial\Omega$. Then the
function $x\mapsto d(x,\partial\Omega):=\inf\{d(x,y) :
y\in\partial\Omega\}$ is a viscosity solution of the equation {\em
(HJ4)}. Moreover, if $M$ is uniformly locally convex and has a
positive injectivity radius, then $d(\cdot,\partial\Omega)$ is the
{\em unique} viscosity solution of this equation.
\end{thm}
\begin{proof}
Let us first check uniqueness. Assume $u,
v:\overline{\Omega}\to\mathbb{R}$ are viscosity solutions of
(HJ4). Since $u$ and $v$ are continuous, and $u=v=0$ on
$\partial\Omega$, we can extend $u$ and $v$ with continuity to the
whole of $M$ by setting $u=0=v$ on $M\setminus\Omega$. It is
enough to see that $u\leq v$ on $\Omega$ (in a similar way, or by
symmetry, $v\leq u$, hence $u=v$). To this end we take any
$\alpha\in (0,1)$ and we check that $\alpha u\leq v$. Indeed,
suppose we had that $\inf\{v(x)-\alpha u(x) :x\in\Omega\}<0$. Pick
$\varepsilon$ with
    $$
    0<2\varepsilon<\min\Big\{\frac{1-\alpha}{2}, \,
    -\inf\{v(x)-\alpha u(x) :x\in\Omega\}\, \Big\}.
    $$
Note that, as $u$ and $v$ are viscosity solutions, we have
$\|\zeta\|_{x}\leq 1$ for every $\zeta\in D^{+}u(x)\cup
D^{+}v(x)$, $x\in\Omega$, so by the mean value Theorem \ref{TVM
Deville} $u$ and $v$ are $1$-Lipschitz. In particular, since
$\Omega$ is bounded we have that $u$ and $v$ are bounded. Then,
according to Proposition \ref{HC prop 1}, there exist $x_{0},
y_{0}\in M$, $\zeta\in D^{+}(\alpha u)(x_{0})$, $\xi\in
D^{-}v(y_{0})$ with
\begin{enumerate}
\item $d(x_{0},y_{0})<\varepsilon$
\item $\|\zeta- L_{x_{0}y_{0}}(\xi)\|_{x_{0}}<\varepsilon$
\item $\inf (v-\alpha u)\geq v(y_{0})-\alpha
u(x_{0})-\varepsilon$.
\end{enumerate}
Taking into account the facts that $u$ and $v$ are $1$-Lipschitz,
and $u=v=0$ on $M\setminus\Omega$, it is easy to see that $(3)$
and the choice of $\varepsilon$ imply that $x_{0},
y_{0}\in\Omega$. Now, since $u$ and $v$ are viscosity solutions we
have that
    $$
    \frac{1}{\alpha}\zeta\in D^{+}u(x_{0}) \implies
    \|\frac{1}{\alpha}\zeta\|_{x_{0}}\leq 1 \implies
    \|\zeta\|_{x_{0}}\leq \alpha,\, \textrm{ and }
    $$
    $$
    \xi\in D^{-}v(y_{0}) \implies \|\xi\|_{y_{0}}\geq 1.
    $$
Now, from $(2)$, and bearing in mind that $L_{x_{0}y_{0}}$ is a
linear isometry, we get that
    $$
    1\leq\|\xi\|_{y_{0}}=\|L_{x_{0}y_{0}}(\xi)\|_{x_{0}}\leq
    \|\zeta\|_{x_{0}}+\varepsilon\leq\alpha+\varepsilon<1,
    $$
a contradiction.

\medskip

Now let us prove that $u:=d(\cdot,\partial\Omega)$ is a viscosity
solution to (HJ4), hence the only one. The property $u=0$ on
$\partial\Omega$ is obvious from the definition, so we only have
to check the conditions on the norms of the vectors of
$D^{-}u(x)$ and $D^{+}u(x)$, for $x\in\Omega$.

\noindent {\bf Step 1.} Take $\xi\in D^{-}u(x)$, $x\in\Omega$. We
have to see that $\|\xi\|_{x}\geq 1$. By Theorem \ref{equivalent
definitions of subdifferential} we can pick a $C^1$ smooth
function $\varphi:M\to\mathbb{R}$ so that $u(y)-\varphi(y)\geq
u(x)-\varphi(x)=0$ for all $y\in M$. Fix $0<\varepsilon<1$. Now,
for every $\alpha$ with $0<\alpha<d(x,\partial\Omega)$, by the
definition of $d(x,\partial\Omega)$ we can take
$x_{\alpha}\in\partial\Omega$ with
    $$
    d(x,\partial\Omega)\geq
    d(x,x_{\alpha})-\frac{\varepsilon\alpha}{4}.
    $$
Next, by making use of Ekeland's approximate Hopf-Rinow type
Theorem \ref{Ekelan's approximated Hopf-Rinow theorem}, we can
find a point $y_{\alpha}\in\Omega$ with
    $$
    d(x_{\alpha},y_{\alpha})<\frac{\varepsilon\alpha}{4}
    $$
and a geodesic $\gamma_{\alpha}:[0,
T_{\alpha}]\to\overline{\Omega}\subset M$ joining
$x=\gamma_{\alpha}(0)$ to
$y_{\alpha}=\gamma_{\alpha}(T_{\alpha})$, and such that
$L(\gamma_{\alpha})=d(x,y_{\alpha})$. Then we have
    $$
    L(\gamma_{\alpha})=d(x,y_{\alpha})\leq
d(x,x_{\alpha})+d(x_{\alpha},y_{\alpha})\leq
d(x,x_{\alpha})+\frac{\varepsilon\alpha}{4}\leq
d(x,\partial\Omega) +\frac{\varepsilon\alpha}{2},
    $$
that is
    $$
    d(x,\partial\Omega)\geq
    L(\gamma_{\alpha})-\frac{\varepsilon\alpha}{2}. \eqno(4)
    $$
Set $v_{\alpha}=d\gamma_{\alpha}(0)/dt\in TM_{x}$, so that
$\gamma_{\alpha}(t)=\exp_{x}(tv_{\alpha})$ and
$\|v_{\alpha}\|_{x}=1$, and define
$z_{\alpha}=\gamma_{\alpha}(\alpha)$. Then we have
\begin{eqnarray*}
& &\varphi(z_{\alpha})-\varphi(x)\leq u(z_{\alpha})-u(x)=
d(z_{\alpha},\partial\Omega)-d(x,\partial\Omega)\leq\\ & &
d(z_{\alpha},\partial\Omega)-L(\gamma_{\alpha})+\frac{\varepsilon\alpha}{2}
\leq d(z_{\alpha},y_{\alpha})+d(y_{\alpha},x_{\alpha})
-L(\gamma_{\alpha})+\frac{\varepsilon\alpha}{2}\leq\\ & &
L({\gamma_{\alpha}}_{|_{[\alpha,
T_{\alpha}]}}))+\frac{\varepsilon\alpha}{2}
-L(\gamma_{\alpha})+\frac{\varepsilon\alpha}{2}=\\ & &
L(\gamma_{\alpha})-\alpha+\frac{\varepsilon\alpha}{2}
-L(\gamma_{\alpha})+\frac{\varepsilon\alpha}{2}=\alpha(-1+\varepsilon),
\end{eqnarray*}
hence
    $$
    \frac{\varphi(z_{\alpha})-\varphi(x)}{\alpha}\leq
    -1+\varepsilon. \eqno(5)
    $$
By the mean value theorem there is $s_{\alpha}\in [0,\alpha]$
such that
    $$
    d\varphi(\gamma_{\alpha}(s_{\alpha}))(\frac{d\gamma_{\alpha}(s_{\alpha})}{dt})=
    \frac{\varphi(z_{\alpha})-\varphi(x)}{\alpha}. \eqno(6)
    $$
By combining$(5)$ and $(6)$, and bearing in mind that
$\|d\gamma_{\alpha}(s)/dt\|_{\gamma_{\alpha}(s)}=1$ for all $s$,
we get that
    $$
    \|d\varphi(\gamma_{\alpha}(s_{\alpha}))\|_{\gamma_{\alpha}(s_{\alpha})}
    \geq 1-\varepsilon \eqno(7)
    $$
for every $\alpha\in (0, u(x))$. Then, since the functions $y\to
d\varphi(y)$ and $(y,\zeta)\to\|\zeta\|_{y}$ are continuous, and
$\gamma_{\alpha}(s_{\alpha})=\exp_{x}(s_{\alpha} v_{\alpha})\to
\exp_{x}(0)=x$ as $\alpha\to 0$, it follows that
    $$
    \|\xi\|_{x}=\|d\varphi(x)\|_{x}=
    \lim_{\alpha\to 0^{+}}
    \|d\varphi(\gamma_{\alpha}(s_{\alpha}))\|_{\gamma_{\alpha}(s_{\alpha})}
    \geq 1-\varepsilon. \eqno(8)
    $$
Finally, by letting $\varepsilon\to 0$ in $(8)$, we deduce that
$\|\xi\|_{x}\geq 1$.

\noindent {\bf Step 2.} Now take $\zeta\in D^{+}u(x)$,
$x\in\Omega$, and let us see that $\|\zeta\|_{x}\leq 1$. This is
much easier. Pick a $C^1$ smooth function $\psi:M\to\mathbb{R}$ so
that $d\psi(x)=\zeta$ and $u(y)-\psi(y)\leq u(x)-\psi(x)=0$ for
all $y\in M$. For each $v\in TM_{x}$ consider the geodesic
$\gamma_{v}(t)=\exp_{x}(tv)$. Since $u=d(\cdot,\partial\Omega)$
is 1-Lipschitz we have that
    $$
    \psi(\gamma_{v}(t))-\psi(\gamma_{v}(0))\geq
    u(\gamma_{v}(t))-u(\gamma_{v}(0))\geq -d(\gamma_{v}(t), \gamma_{v}(0))=-t,
    $$
hence
    $$
    \frac{\psi(\gamma_{v}(t))-\psi(\gamma_{v}(0))}{t}\geq -1
    $$
for all $t>0$ small enough, and
    $$
    d\psi(x)(v)=\lim_{t\to 0^{+}}\frac{\psi(\gamma_{v}(t))-
    \psi(\gamma_{v}(0))}{t}\geq -1. \eqno(9)
    $$
As $(9)$ holds for every $v\in TM_{x}$, we conclude that
$\|\zeta\|_{x}=\|d\psi(x)\|_{x}\leq 1$.
\end{proof}


\medskip

\begin{center}
{\bf {\small ACKNOWLEDGEMENTS}}
\end{center}
\noindent We thank Luis Guijarro for several valuable discussions
concerning the results of this paper.

\medskip


\end{document}